\documentclass{article}
\usepackage{amssymb}
\usepackage{amsmath}
\usepackage{amsfonts}
\usepackage{graphicx}

\newtheorem{theorem}{Theorem}
\newtheorem{definition}[theorem]{Definition}

\newtheorem{corollary}[theorem]{Corollary}
\newtheorem{lemma}[theorem]{Lemma}
\newtheorem{proposition}[theorem]{Proposition}
\newtheorem{remark}[theorem]{Remark}

\newcommand{\fin}{\hfill\mbox{$\quad{}_{\Box}$}}
\newcommand{\fineq}{\vspace{-.75cm$\fin$}\par\bigskip}

\begin{document}

\title{Electron beams: partially flat solutions of a nonlinear elliptic
equation with a singular absorption term}
\author{Jes\'{u}s Ildefonso D\'{\i}az \\
Instituto de Matem\'{a}tica Interdisciplinar, \\
Universidad Complutense de Madrid, \\ Plaza de las Ciencias, n%
${{}^o}$
3, Madrid, 28040, Spain \\ jidiaz@ucm.es}
\maketitle

\begin{center}
\textit{Dedicated to my Master, Ha\"{\i}m Brezis, who trusted me and gave me
wings to fly higher.}
\end{center}

\textit{\ }

\textbf{Abstract. }In the so-called Child-Langmuir law, established since 1911,
an electron beam is formed linking two electrodes, which are assumed to be
two parallel plates of area $A$, separated to a finite distance $D.$ When $%
D\ll \sqrt{A},$ \textquotedblleft edge effects\textquotedblright\ are
negligible and the modelling is reduced to a nonlinear boundary problem for a
singular ordinary differential equation\ in which a constant coefficient
(the generated electric current  $j$) must be found in order to get
simultaneously Dirichlet and Neumann homogeneous boundary conditions in one
of the extremes. If $D>\sqrt{A},$ then the problem becomes much more
difficult since the \textquotedblleft edge effects\textquotedblright\ arise
in the plane $(x,y)$ and the electric current (now $j(x)$ due to the
presence of a very large perpendicular magnetic field) must be determined in
order to get solutions $u(x,y)$ of a singular semilinear equation which are
partially flat ($u=\frac{\partial u}{\partial n}=0$ on a part of the
boundary). In this paper, we offer a rigorous mathematical treatment of some
former studies (Joel Lebowitz and Alexander Rokhenko (2003) and Alexander
Rokhenko (2006)), where several open questions were left open: for instance,
the need for a singularity of $j(x)$ near the cathode edge to get such
partially flat solutions.

\bigskip 



\section{Introduction}

In 1911, the American physicist Clement D. Child (1868 -- 1933) proposed in 
\cite{Child} an experimental law (\textquotedblleft Child's
law\textquotedblright ) that he deduced after modelling the electric current
that flows between the plates of a vacuum tube: the main components in
electronics from about 1905 to 1960, when transistors and integrated
circuits mostly supplanted them. For some time (1907-1908), Child was visiting 
J.J. Thompson (1856-1940), 1906 Nobel Prize in Physics, who built the cathode ray
tubes,  proving in 1897 the existence and charge of the electron. Child's law is
still a staple of textbooks treating charged particle motion in vacuum and
in solids. After normalizing the variable and the unknown, and by
considering the case in which the anode and cathode are simulated by planes at a
finite distance, Child proposed a singular nonlinear ordinary differential equation
 which perhaps is one of the older singular equations considered in
the literature. Its formulation was the following: find $j>0$ (modelling the
maximum generated electric current) to get a positive solution of
the nonlinear boundary problem 
\begin{equation}
\left\{ 
\begin{array}{lcc}
-u^{\prime \prime }(y)+\displaystyle\frac{j}{\sqrt{u(y)}}=0 & \quad y\in
(0,1), &  \\ 
u(0)=0\quad u(1)=1, &  & 
\end{array}%
\right. 
\end{equation}%
such that 
\begin{equation*}
u^{\prime }(0)=0\text{ and }u>0\text{ on }(0,1).
\end{equation*}%
Notice that since $u(0)=0$ the singularity of the equation is present on one
of the boundaries ($y=0$), representing the cathode. Condition $u^{\prime
}(0)=0$  says that the electron flux is assumed to start from rest,
since the electric field is zero there. This law was later considered in a
series of papers (see, e.g., \cite{Langmuir}), and improved by Irving
Langmuir (1881-1957), 1932 Chemistry Nobel Prize. Since then, the above law
is called the \textquotedblleft Child-Langmuir law\textquotedblright . In
physical terms, it is known as the \textquotedblleft three-halves-power
law\textquotedblright\ (see Remark \ref{Rmk Child law}), since it says that
the electric current $J$ is proportional to the $3/2$-power of the voltage, $%
V^{3/2}$ (in contrast to Ohm's law on electric circuits for which $J=kV$
for some $k>0$). We will see, in Section 3, that in the renormalized
formulation

\begin{equation*}
j\equiv \frac{4}{9}\text{ and }u(y)=y^{4/3}\text{ for }x\in (0,1),
\end{equation*}%
in contrast with the 1823 Ohm's law [$u(y)=ky$, for some $k>0$].

In the Child-Langmuir modelling, the electrodes are assumed to be two
parallel plates of area $A$, separated by a finite distance $D$ ($D=1$ after
some change of variables). They assumed that $D\ll \sqrt{A},$ so that
\textquotedblleft edge effects \textquotedblright are negligible. 

The study of the problem in which edge effects are taken into account is
much more difficult and was a subject of research until our days (see, e.g.,
the survey \cite{Zhang}). A more realistic modelling, already considered by
many authors (see, e.g., \cite{Rokhenko-Lebowitz2003}, \cite%
{Rokhenko-Lebowitz2005}), takes place in the plane $(X,y)$, in which the cathode is only
an interval $(-a,a)\times \{0\}$ of the device (the set of points $(X,y)\in (-a-b,a+b)\times
(0,1)$, for a given $b>0$), and it is assumed the presence of a strong magnetic
field. This makes the electron beam well confined in such a way that the
electric current is only $X$-dependent, $J(X,y)=j(X)$. The primitive
question regarded by the Child-Langmuir law can now be stated in the
following terms as an overdetermined problem: given $a,b>0,$ find sufficient
conditions on a $x$-dependent function $j:(-a-b,a+b)\mathbb{\rightarrow
\lbrack }0,+\infty )$, with 
\begin{equation}
\left\{ 
\begin{array}{lr}
j(X)>0 & \text{if }X\in (-a,a),\text{ }j\in L_{loc}^{1}(-a,a), \\ [.15cm]
j(X)=0 & \text{if }X\in (-a-b,-a)\cup (a,a+b),%
\end{array}%
\right.   \label{Hypo j for a}
\end{equation}%
to get the solvability of the singular nonlinear boundary value
problem 
\begin{equation}
\widehat{P}_{a,b,j}=\left\{ 
\begin{array}{lr}
-\Delta u+\frac{j(X)}{\sqrt{u}}=0 & X\in (-a-b,a+b),\text{ }y\in (0,1), \\ [.15cm]
u(X,0)=0 & X\in (-a-b,a+b), \\ [.15cm]
u(X,1)=1 & X\in (-a-b,a+b), \\ [.15cm]
u(\pm (a+b),y)=y & y\in (0,1),%
\end{array}%
\right.   \label{BVP}
\end{equation}%
with the additional conditions 
\begin{equation}
\widehat{AC}_{a,b}=\left\{ 
\begin{array}{lr}
\frac{\partial u}{\partial y}(X,0)=0 & X\in (-a,a)\text{, } \\ [.15cm]
u(X,y)>0 & X\in (-a-b,a+b),\text{ }y\in (0,1).%
\end{array}%
\right.   \label{AC}
\end{equation}

The more important fact considered in the literature is the study of the
adaptation of profiles $u(\cdot ,y)$ from the external profile (which is
assumed to be given by  $u(b,y)=y$, as in Ohm's law) to the profile in the
center of the cathode $u(0,y),$ where the authors expect to have the profile
correspondent to the Child-Langmuir law, $u(0,y)=y^{4/3}$ for $y\in (0,1).$

A very delicate question is the behaviour of $j(X)$ near $X=\pm a$ since it
must generate complicated edge effects (for instance, by the \textit{strong
maximum principle}, we know that $\frac{\partial u}{\partial y}(X,0)>0$ if $%
(-a-b,-a)\cup (a,a+b)$, since $u$ is harmonic on $\left( (-a-b,-a)\cup
(a,a+b)\right) \times (0,1)$).

To simplify the formulation, due to the symmetry of the problem, it is usual
to consider only (see, e.g., \cite{Rokhenko-Lebowitz2003}, \cite%
{Rokhenko-Lebowitz2005}, \cite{Rokhlenko}) the transmission of the profiles $%
u(X_{0},y),$ $y\in (0,1)$, when $X_{0}\in \lbrack 0,a+b]$ in one of the two symmetrical parts of the device (for example, the one on the right). To give 
the edge behaviour (on the boundary of the cathode) a central presence in
the formulation, it is traditional to introduce the change of variable $x=X-a
$, $X\in (0,a+b)$ and then $x\in (-a,b)$ for a given $b>0$. In this way, $(-a,0)$ becomes the center of the cathode and $(0,0)$ the edge of the right side of the cathode. We arrive then
to the formulation which will be considered in this paper: given $a,b>0$, we
consider the problem of finding the coefficient $j(x)$ such that there
exists a solution of the singular boundary value problem on $\Omega
=(-a,b)\times (0,1)$ 
\begin{equation}
P_{a,b,j}=\left\{ 
\begin{array}{lr}
-\Delta u+\frac{j(x)}{\sqrt{u}}=0 & \text{in }\Omega , \\  [.15cm]
u(x,0)=0 & x\in (-a,b), \\ [.15cm]
u(x,1)=1 & x\in (-a,b), \\ [.15cm]
u(-a,y)=y^{4/3} & y\in (0,1), \\ [.15cm]
u(b,y)=y & y\in (0,1),%
\end{array}%
\right.   \label{Pa b j}
\end{equation}%
satisfying the additional conditions
\begin{equation}
AC_{a,b}=\left\{ 
\begin{array}{lr}
\frac{\partial u}{\partial y}(x,0)=0 & x\in (-a,0)\text{, } \\ [.15cm]
u>0 & \text{in }\Omega .%
\end{array}%
\right.   \label{AC ab}
\end{equation}

Notice that the additional conditions imply a failure of the
\textquotedblleft unique continuation property\textquotedblright\ and the
Hopf-Oleinik lemma related to the strong maximum principle (see, e.g., the exposition made in \cite
{DiazHernandezBeyond}): indeed, $u>0$ in $(-a,b)\times (0,1)$ although $u=%
\frac{\partial u}{\partial y}=0$ on $(-a,0)\times \{0\}$. The additional
condition, $\frac{\partial u}{\partial y}(x,0)=0$ for $x\in (-a,0)$,
represents the vanishing of the electric field on the half of the cathode under consideration.   This type of behaviour arises in many \textit{free boundary problems }(\cite{D}) for semilinear equations with a non-Lipschitz absorption. It is
quite similar to the requirement of the so-called \textit{flat solutions }%
arising, for instance, in the study of linear (\cite{Diazambiguous}, \cite%
{DiazHernandezBeyond}) and nonlinear (\cite{DHYliasov}) Schr\"{o}dinger
equations. Here, the constraint on the function $j$ and the
singularity of the nonlinear absorption term ($1/\sqrt{u}$) generate some
important additional difficulties.

\begin{figure}[htp]
\begin{center}
\includegraphics[width=7cm]{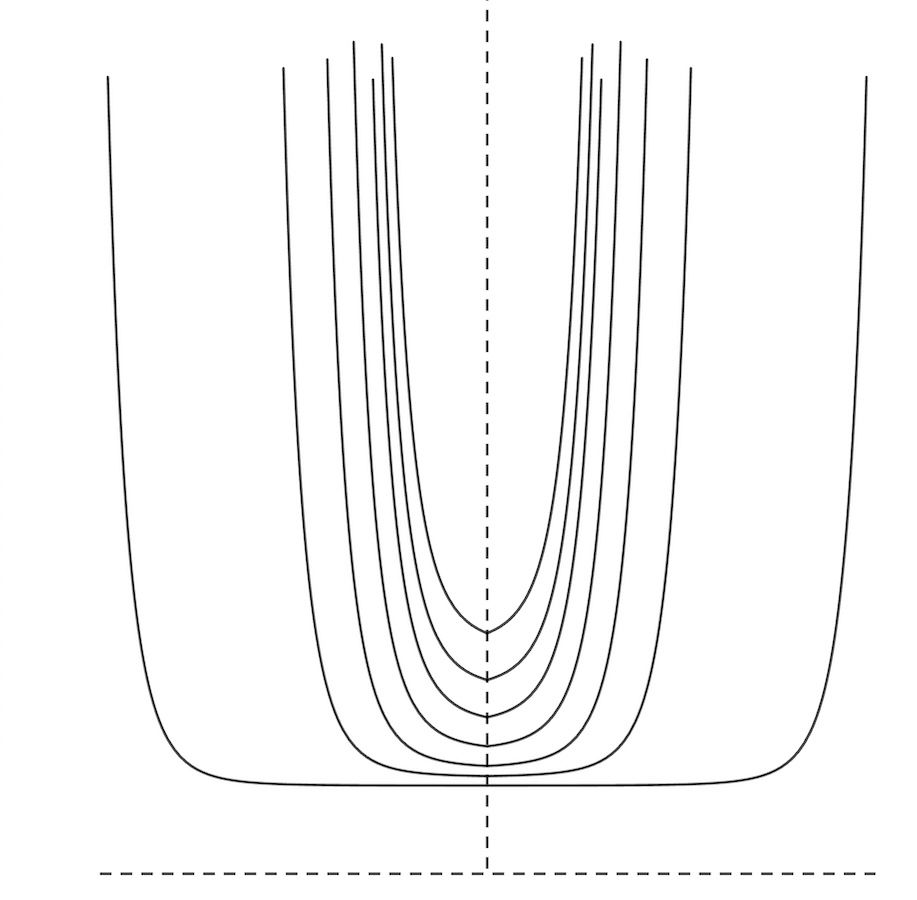}\\ 
\caption{Wings of $j(x)$ in some
numerical experiments for different values of the cathode width $a$: adapted from \protect\cite{Umstattd}.}
\label{Fig alas}
\end{center}
\end{figure}

Despite the relevance of the space-charge-limited flows in many
applications (vacuum and solid electronic devices, electron guns, particle
accelerators, high current diodes, tubes, thyristors, etc., see, e.g., \cite%
{Kirsten}, \cite{H2002}, \cite{Hoffman}) the study of 2-d simple geometries
as here considered presents important difficulties, and it was the
benchmark for almost a century, in particular for the study of the generated current $j(x)$ in the vicinity of the cathode edges (see Figure \ref{Fig alas}). Many previous studies, based on
numerical and asymptotic methods (see, e.g., \cite{Watrous}, \cite{Rokhenko-Lebowitz2003}, \cite{Rokhenko-Lebowitz2005}, \cite%
{Rokhlenko} and their many references) claim that due to the important
constraint on $j(x)$ (independence with respect to $y$) it is needed to assume that%
\begin{equation*}
\lim_{x\nearrow 0}j(x)=+\infty ,
\end{equation*}%
(see, in Figure \ref{Fig alas}, the singularities of the function $j(x)$ at the extremes $\pm a$, for different widths $a$ of the cathode). A rigorous mathematical proof of it was left
as an open problem in the cited references.


One of the main goals of this paper is to prove, and make it precise, a
conjecture proposed in many numerical analysis experiences (see, e.g., \cite%
{Rokhlenko} and its references) concerning the formation of
\textquotedblleft wings\textquotedblright\ in the electric current (i.e., the
formation of a singularity of $j(x)$ near the cathode edge). In addition, we will
prove a stronger version of the positivity of the solution: in fact, the solution is non-degenerate
(see subsection 4.4).

\begin{theorem}
\textbf{\label{TheomMain}}\textit{There exists }$A_{0},b_{0}>0$\textit{\ and 
}$\beta _{0}\in (0,\frac{1}{2}),$\textit{\ such that, if }$b\geq b_{0}>0,$ 
\textit{and we assume} 
\begin{equation}
\left \{ 
\begin{array}{lr}
j(x)= \frac{A}{(-x)^{\beta }}, & \text{\textit{for} }\mathit{%
x\in (-a,0),} \\[0.15cm]
\mathit{j(x)=0} & \text{\textit{for} }x\in (0,b),%
\end{array}%
\right.   \label{Hypo j supersolution}
\end{equation}%
\textit{with}%
\begin{equation}
0\leq \beta <\beta _{0}\text{ }\mathit{and}\text{ }A\in (0,A_{0})\mathit{.}
\end{equation}%
\textit{\ Then there exists a weak solution} $u\in L^{2}(\Omega ;\delta )$ with%
\textit{\  \ }$\delta =d((x,y),\partial \Omega ),$ \textit{of problem }$%
P_{a,b,j}$ and the additional conditions $AC_{a,b}$. Moreover\textit{\ } 
\begin{equation*}
0<\underline{C}\delta (x,y)^{4/3}\leq u(x,y)\leq \overline{C}%
\delta (x,y)^{\overline{\alpha }}\text{ a.e. }(x,y)\in (-a,0)\times (0,1),
\end{equation*}%
\textit{and} $Cy^{\alpha }\leq u(0,y),$ $y\in (0,1),$ with $%
\alpha =\frac{2}{3}(2-\beta )$, $1<\overline{\alpha }\leq \frac{4}{3}$, $C$,$\underline{C}$,$\overline{C}>0$.  Finally, $u$ is unique in the class of non-degenerate solutions.$\fin$
\end{theorem}

\bigskip

The proof will be made by constructing suitable super and subsolutions for
several auxiliary problems that are obtained by matching some strategic functions on
different subdomains. In particular, we will need to consider the nonlinear
singular eigenvalue type problem 
\begin{equation}
\left\{ 
\begin{array}{lr}
-U^{\prime \prime }(s)+\frac{C}{\sqrt{U(s)}}=\lambda U(s) & s\in (-R,R), \\ [.2cm]
U(\pm R)=0, & 
\end{array}%
\right.   \label{Problem Psi}
\end{equation}

\begin{figure}[htp]
\begin{center}
\includegraphics[width=9cm]{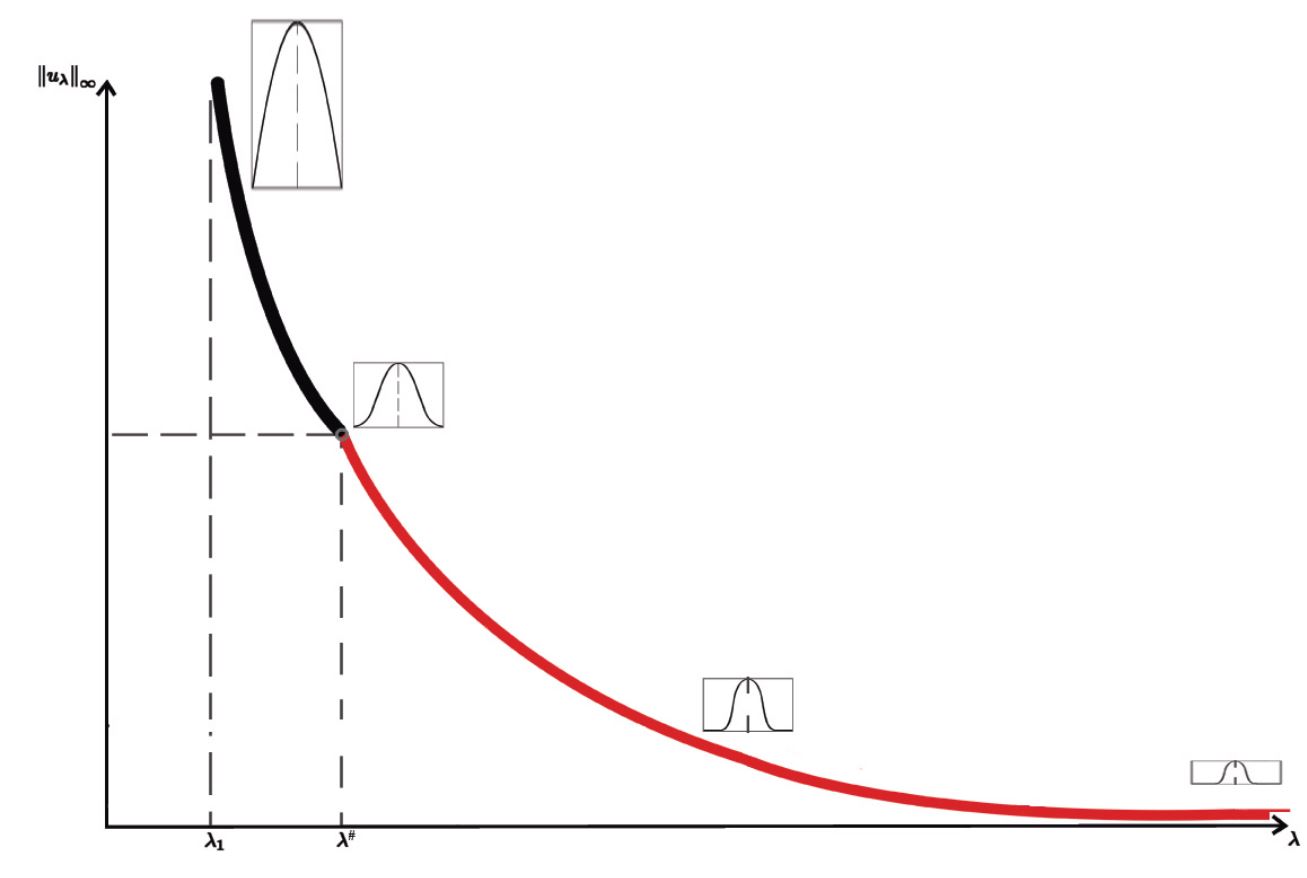}\\ 
\caption{Bifurcation from infinity and critical values of $\protect\lambda$ for a flat solution.}
\label{Fig Canarias}
\end{center}
\end{figure}
\noindent for which the global bifurcation diagram (in terms of the
parameter $\lambda $) will be completely characterized (see Figure \ref{Fig Canarias}). Here the positive
constants $C$ and $R$ are given. We will show: i) there is a bifurcation
from the infinity for $\lambda $ near $\lambda _{1}(R)$ (the first
eigenvalue of the linear problem with $C=0$), ii) the bifurcation curve is
strictly decreasing (which implies the uniqueness of the nonnegative solution $U
$) and iii) the curve is not $C^{1}$ for a suitable value $\lambda =\lambda
^{\ast }>\lambda _{1}(R)$ corresponding to a \textquotedblleft flat
solution\textquotedblright\ (i.e. the solution $U$ is such that $U^{\prime
}(\pm R)=0$ and $U(s)>0$ for any $s\in (-R,R)$). This extends several
results in the previous literature (see, e.g., \cite{Di-Hern-Mance}, \cite%
{DHPortugalia}, and their references).

\begin{figure}[htp]
\begin{center}
\includegraphics[width=9cm]{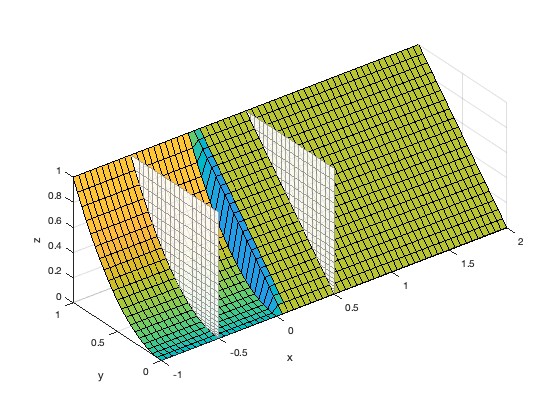}\\ 
\caption{Progressive changes of profiles $u(\cdot,y).$}
\label{Fig Profiles}
\end{center}
\end{figure}

Theorem \ref{TheomMain} confirms many numerical experiences that lacked a rigorous proof. The progressive changes in the profiles can be illustrated
by Figure \ref{Fig Profiles}.

We point out that the nature of the problem changes radically when the
singular nonlinearity arises on the other side of the equality (as a forcing
term): see, e.g. \cite{CrandallRT} and \cite{DHRako}. Most of the
motivations to consider singular absorption terms have their foundation in fields such as chemical
reactions, non-Newtonian flows and other applied subjects very different from
the one of space-charge-limited flows (see, e.g., \cite{Fulks}, \cite{DMO87}%
, \cite{Choi-Lazer-Mckenna} and \cite{D}). A survey on singular elliptic
equations can be found in \cite{Hernandez-Mancebo}.

We will need to construct several super and subsolutions satisfying suitable
qualitative properties. This will require matching different functions. We
point out that two remarkable matching studies in the literature were those of
G. Gamov in \cite{Gamow} (when proving the tunneling effect) and the one\
by H. Brezis in \cite{Brezis-soporte} (proving the compact support of
solutions of suitable variational inequalities: see the exposition made in 
\cite{DGaeta}). Here we will need to use a matching method that is not the usual $H^{1}$-matching
criterion but a looser one that allows the formation of singularities as long as they have a certain sign (see Remark \ref{Rmk Matching picos} below).

The organization of this paper is the following: Section 2 is devoted to the
modelling of the problem. In particular, we will make mention to some
optimization criterion on $j(x)$ (see Remark \ref{Rm Optimization}). The
one-dimensional case, including the Child-Langmuir law, is revisited in
Section 3. The elegant study of the one-dimensional case is due to H. Brezis
(personal communication in 2004). We also consider the non-autonomous case $%
j=j(x)$ since it will be applied later. The main goal of this work is the study of the two-dimensional formulation that will be developed in the last Section 4 divided into different subsections. We will give details on the
construction of a positive subsolution and a supersolution, both satisfying the
additional condition, and prove the uniqueness of solutions in the class of
non-degenerate solutions.

\section{On the modelling of the 2-d space-charge-limited flow}

\noindent Following \cite{Rokhenko-Lebowitz2003} (see also \cite{Pollard}),
a general formulation of the problem under consideration could start with the
consideration of the stationary Maxwell system of equations for the electric
and magnetic fields $(\mathbf{E,B)}$ defined on a set $\widetilde{\Omega }%
\subset \mathbb{R}^{3},$ $\widetilde{\Omega }=\mathbb{R\times }(0,D)\times $ 
$\mathbb{R}$, with $D>0$ given, separating two conducting electrodes placed
on the planes $Y=0$ (cathode) and $Y=D$ (anode) [with $\partial \widetilde{%
\Omega }=\Gamma _{0}\cup \Gamma _{1}$ ]:

\begin{equation}
\left\{ 
\begin{array}{lr}
\nabla \cdot \mathbf{E}=\displaystyle\frac{\rho }{\epsilon _{0}} &  \\ [.2cm]
\nabla \cdot \mathbf{B}=0 &  \\ [.2cm]
\nabla \times \mathbf{E}=\mathbf{0} &  \\ [.2cm]
\nabla \times \mathbf{B}=\mu _{0}\mathbf{J.} & 
\end{array}%
\right.  \label{Maxwell equations}
\end{equation}%
Here $\rho (X,Y,Z)$ is the electron-charge density, $\epsilon _{0}$ the
free space permittivity and $\mathbf{J}(X,Y,Z)$ denotes the current density
(here we are using a variable notation $(X,Y,Z)$ different from the one of  
the Introduction). Since we are assuming that $\rho $ is stationary we
get that 
\begin{equation*}
div~\mathbf{J}=0\text{ in }\widetilde{\Omega }.
\end{equation*}

Now we assume that the \textit{cathode} is in the $(X,Z)$ plane, $Y=0$, and
it has a width $2A$. We also assume that there is a \textit{very strong
magnetic field} $\mathbf{B}$, which is perpendicular to the electrodes ($%
\mathbf{B}$ is in the $Y$-direction), inhibiting the transversal components
of the electron velocities $v(X,Y)$, and then $\rho (X,Y,Z)=\rho (X,Y)\chi
_{\{\left\vert X\right\vert \leq A\}}(X,Y)$ , where $\chi _{\{\left\vert
X\right\vert \leq A\}}(X,Y)$ is the characteristic function 
\begin{equation*}
\chi _{\{\left\vert X\right\vert \leq A\}}(X,Y)=\left\{ 
\begin{array}{cc}
1 & \text{if }\left\vert X\right\vert \leq A \\ 
0 & \text{otherwise.}%
\end{array}%
\right. 
\end{equation*}%
The fact that $\rho $ is not constant is the reason to call this type of
problems as \textit{space charge}. Due to the assumption on $\mathbf{B}$, we
know that \textit{the potential} $U$ of the \textit{electric field}\textbf{\
(}$\mathbf{E=-\nabla }U\mathbf{)}$ is $Z$-independent, i.e., $U=U(X,Y)$.

We assume that \textit{the emitted electrons leave the cathode with zero
velocity}\textbf{\ }and thus, if we take $U=U(X,Y)=0$ in the cathode, the
total \textit{mechanical energy} is $E_{0}=0$. If $e$ and $m$ represent the 
\textit{charge} and the \textit{mass} of the electron, the \textit{%
conservation of the mechanical energy} leads to the equation%
\begin{equation}
\frac{m}{2}v^{2}(X,Y)=eU(X,Y).  \label{Conservation}
\end{equation}%
i.e.,

\begin{equation*}
v(X,Y)=\sqrt{\frac{2eU(X,Y)}{m}}.
\end{equation*}%
Remember that the mechanical force (by a negative charge) is given by $%
\mathbf{F}=(-e)\mathbf{E=-}e(-\nabla U)=e\nabla U$ and thus the potential
energy is $-eU.$ From this, we deduce that the current density is only
dependent on $X$, $\mathbf{J}(X,Y,Z)=J(X)\chi _{\{\left\vert X\right\vert
\leq A\}}(X,Y)\mathbf{e}_{2}$ and determines the velocity of electrons. We
also recall that 
\begin{equation*}
\displaystyle \mathbf{J}=-\rho v\mathbf{e}_{2}:=J(X)\chi _{\{\left\vert X\right\vert \leq
A\}}(X,Y)\mathbf{e}_{2}.
\end{equation*}%
Then%
\begin{equation*}
\rho (X,Y)=-\frac{J(X)\chi _{\{\left\vert X\right\vert \leq A\}}(X,Y)}{\sqrt{%
\frac{2eU(X,Y)}{m}}}.
\end{equation*}

We introduce now the \textit{dimensionless variables }(once again, the
notation is different from the one used in the Introduction)%
\begin{equation*}
x=\frac{X}{D},y=\frac{Y}{D}\text{ and }a=\frac{A}{D},
\end{equation*}%
and the \textit{dimensionless functions}%
\begin{equation*}
u(x,y)=\frac{U(X,Y)}{V}\text{ and  }j(x)=\frac{9}{4}\sqrt{\frac{m}{2e}}\frac{%
J(X)}{\epsilon _{0}V^{3/2}}.\text{ }
\end{equation*}%
Then, the first equation of the Maxwell system and the conservation of the
mechanical energy lead to the singular nonlinear Poisson equation%
\begin{equation*}
\Delta u(x,y)=-4\pi \rho (x,y)=\frac{j(x)}{\sqrt{u(x,y)}},
\end{equation*}%
with $\rho (x,y)=0$ (i.e. $j(x)=0$) if $\left\vert x\right\vert >a$.

\bigskip

\begin{remark}\rm
\label{Rmk Child law}In the dimensionless process, it is possible to choose a
different new expression for $J(X).$ So, for instance, in \cite%
{Rokhenko-Lebowitz2003} it was used the function $\widetilde{j}(x)=\frac{4}{9%
}j(x),$ so that the one-dimensional case corresponds to $\widetilde{j}%
(x)\equiv 1.$ In any case, the \textquotedblleft three-halves-power
law\textquotedblright\ Child-Langmuir law says that in the one-dimensional
case%
\begin{equation}
J=\frac{4}{9}\epsilon _{0}\sqrt{\frac{2e}{m}}\frac{V^{3/2}}{D^{2}}.
\label{Child-Langmuir law}
\end{equation}
\end{remark}

\begin{remark}\rm
Most of the results of the following sections of this paper can be applied
in the framework of two concentric cylindrical diodes (see the formulation
considered in \cite{Chen-Dickens-Choi-Kristiansen}). A different modelling
was considered in \cite{Budd-Friedman} where they study a section through a
thin wire at high potential contained within an earthed conductor at zero
potential.
\end{remark}

\begin{remark}\rm
The study of the mathematical treatment of this problem was suggested to me
by Ha\"{\i}m Brezis, during my visit to his Department of the University of
Paris VI, in \ March 2003. Brezis attended a seminar by Joel
Lebowitz (on his article \cite{Rokhenko-Lebowitz2003}), at his workplace,
Rutgers University (where Brezis would later have a contract as
Distinguished Visiting Professor from 2004 until the date of his death in
2024). I received a more concrete formulation in the fax of Brezis on May,
21, 2004 (containing the Theorem \ref{thm Brezis} below) and I sent him a
first set of my results on the two-dimensional problem in June and September of 2004. \ We made several working sessions during
my stay in Paris, as Invited Professor of Paris VI, from February 15, to
March 14 of 2005, but without finding a positive subsolution (see subsection
4.2 below). In this occasion, we received a visit from Joel Lebowitz to
Paris VI University, and we exchanged opinions in a joint session 
the three of us together (see Remark \ref{Rm General formulation} below and the picture at the end of this paper). I
produced a new draft, sent to both of them, on June 27, 2006. I made a
public presentation of the results in the Opening Workshop of the European
FIRST Project on June 30, 2009, in Orsay (University of Paris
Sud). Later, I presented a first version of the results of this
article at the conference, \textquotedblleft Recent and new perspectives in
Nonlinear Analysis\textquotedblright , held in Urbino (Italy), November 3-4,
2022 and at the IMEIO PhD Course, of the Universidad Polit\'{e}cnica de
Madrid, on January 16, 2023. 
\end{remark}

\begin{remark}\rm
\label{Rm General formulation}(Suggestion by Joel Lebowitz, Paris, 2005).
For other different geometrical electrode shapes, we can assume that the open
(possibly unbounded) set $\widetilde{\Omega }\subset \mathbb{R}^{3}$ has
two components defining its finite boundary $\partial \widetilde{\Omega }%
=\Gamma _{0}\cup \Gamma _{1}$ and the problem is to find a pair $(\mathbf{J}%
(x,y,z),u(x,y,z))$ satisfying (in some weak sense) the following set of
conditions:%
\begin{equation*}
div\mathbf{J}=0\text{ in }\widetilde{\Omega },\text{ }  
\end{equation*}%
\begin{equation*}
u(x,y,z)\geq 0\text{ in }\widetilde{\Omega },
\end{equation*}%
\begin{equation*}
\sqrt{u}\Delta u=j(x,y,z)\text{ in }\widetilde{\Omega },
\end{equation*}%
\begin{equation*}
j(x,y,z)=\left\vert \mathbf{J}(x,y,z)\right\vert ,  
\end{equation*}%
\begin{equation*}
\frac{\mathbf{J}(x,y,z)}{\left\vert \mathbf{J}(x,y,z)\right\vert }=\frac{%
\nabla u(x,y,z)}{\left\vert \nabla u(x,y,z)\right\vert }  
\end{equation*}%
\begin{equation*}
u\left\vert _{\Gamma _{0}}\right. =0\text{ and }u\left\vert _{\Gamma
1}\right. =1.
\end{equation*}
\end{remark}

\begin{remark}\rm
\label{Rm Optimization}As a matter of fact, already in the pioneering works
by Child and Langmuir it was mentioned that the interesting case corresponds
to a current density $J^{\ast }$ defined as the \textquotedblleft largest
current\textquotedblright\ that can be emitted without time-dependent
behavior and thus independent of any thermodynamic effect on the cathode.
This will be illustrated later for the one-dimensional case (when justifying
why $j=4/9$). For a general formulation, as the one presented in Remark \ref%
{Rm General formulation}, this kind of optimization criterion can be stated
as follows (suggestion by Joel Lebowitz, Paris, 2005): from the first of the above conditions
we deduce that there exists a constant $\alpha =\alpha (\mathbf{J)}
$ such that for any simple curve $\Gamma $ (when $\widetilde{\Omega }$ is
unbounded $\Gamma $ can be unbounded without multiple intersection points)
we have that%
\begin{equation*}
\int_{\Gamma }\mathbf{J}\text{\textperiodcentered }\mathbf{n}d\sigma =\alpha
.
\end{equation*}%
We define\textbf{\ }the set $C$ of \textquotedblleft admissible
solutions\textquotedblright\ $(\mathbf{J}(x,y,z),u(x,y,z))$ by means of the
set of conditions given in Remark \ref{Rm General formulation} and then the
problem is now to find $(\mathbf{J}^{\ast },u^{\ast })\in \mathcal{C}$ such
that%
\begin{equation*}
\alpha (\mathbf{J}^{\ast })=\max_{(\mathbf{J},u)\in \mathcal{C}}%
\alpha (\mathbf{J})\text{.}
\end{equation*}
\end{remark}

\begin{remark}\rm
\label{Rm Onedimensional Open problem}As already mentioned in the
Introduction, when $a=+\infty $ the problem can be associated with a
one-dimensional formulation for which we get $\mathbf{J}^{\ast }(x,y,z)=%
\frac{4}{9}\mathbf{e}_{2}$ and, as we will prove in the next Section, $%
u^{\ast }(y)=y^{4/3}.$ An interesting open question can be stated in this
framework\textbf{:} is it true that if $(\mathbf{J},u)\in \mathcal{C}$ then
necessarily $\mathbf{J}(x,y,z)=j\mathbf{e}_{2}$ with $j$ constant and $u=u(y)
$?
\end{remark}

\begin{remark}\rm
\label{Rm Boundary condition at the infinity}In the case of
space-charge-limited flows $(0<a+b<+\infty )$ the boundary condition $u(\pm
(a+b),y)=y$ for $y\in (0,1)$ was first proposed in \cite%
{Rokhenko-Lebowitz2003}, since the external electric field $\mathbf{E}$ must
behave, at least at very long distances, in a similar way to the case
without any cathode ($a=0$).
\end{remark}

\section{The one-dimensional Child-Langmuir law revisited}

\subsection{The autonomous case}

The problem associated with the one-dimensional Child-Langmuir law can be
formulated in the following terms: find $j>0$ in
order to get a function $u\in W^{2,1}(0,1)$ solving the boundary problem 
\begin{equation}
\left\{ 
\begin{array}{lcc}
-u^{\prime \prime }(y)+\displaystyle\frac{j}{\sqrt{u(y)}}=0 & \quad y\in
(0,1),\text{ } &  \\ [.3cm]
u(0)=0\quad u(1)=1, &  & 
\end{array}%
\right.   \label{e1}
\end{equation}%
and such that

\begin{equation}
u^{\prime }(0)=0\text{ and }u>0.  \label{flat condition}
\end{equation}%
We recall that, in this physical case, the only interest is in positive
solutions $u(y)>0$ for any $y\in (0,1)$ but the problem also arises in other
fields (for instance, in Chemical Engineering) where the possibility of
having \textit{a dead core} (a subinterval $(0,\xi )\subset (0,1)$ on which $u=0
$) makes perfect sense for suitable values of $j$.

\noindent Here $W^{2,1}(0,1)$ denotes the usual Sobolev space requiring that 
$u^{\prime \prime }\in L^{1}(0,1)$ (see, e.g., \cite{B2011}). Since $%
W^{2,1}(0,1)\subset C^{1}([0,1])$ condition (\ref{flat condition}) is well
justified

We point out that since the singular nonlinear term $\frac{1}{\sqrt{u}}$ is
neither monotonically non-decreasing, nor Lipschitz continuous, even in
the case in which $j>0$ is given, we cannot apply general results in the
literature on the uniqueness of solutions. 

As a matter of fact, given $j>0,$ some different notions of solution of (%
\ref{e1}) can be introduced (being compatible with the supplementary
condition (\ref{flat condition})). A natural notion of solution arises when
problem \eqref{e1} is understood in the framework of the Calculus of
Variations. We define 
\begin{equation*}
K=\left \{ {u\in H^{1}(0,1)\text{ such that }u(0)=0,u(1)=1\text{ and }u\geq 0%
\text{ on }(0,1)}\right \} ,
\end{equation*}%
and

\begin{equation*}
J(u)=\int_{0}^{1}\left( \frac{1}{2}\ u^{\prime }(y) ^{2}+2j\sqrt{%
u(y)}\right) dy,
\end{equation*}%
and then we call \textit{variational solution of }(\ref{e1}) to a function $%
u\in K$ such that

\begin{equation}
J(u)=\min_{v\in K}\left\{ \displaystyle\int_{0}^{1}\left( \frac{1}{2}
v^{\prime }(y)^{2}+2j\sqrt{v(y)}\right) dy\right\} .  \label{e2}
\end{equation}%
As said before, one of the reasons why the singular problem is quite
relevant in the applications is for the possible occurrence of a \textit{%
free boundary} (according to the values of $j>0$) associated to a solution $u$
of (\ref{e1}). In the one-dimensional setting, the \textit{free boundary} is
given simply by a point $\xi \in \lbrack 0,1)$ such that 
\begin{equation}
\left\{ 
\begin{array}{llc}
u(y)=0\quad \text{if }y\in \lbrack 0,\xi ],\quad  &  &  \\ [.15cm]
u(y)>0\quad \text{if }y\in (\xi ,1],\text{ and }u^{\prime }(\xi )=0. &  & 
\end{array}%
\right.   \label{free boundary condition}
\end{equation}%
Notice that the supplementary condition (\ref{flat condition}) corresponds
to the case $\xi =0$. Many authors used to say that, in that case, $u$ is a
\textquotedblleft \textit{flat solution}\textquotedblright\ (on $y=0$). This
is the more interesting case for us, since the electric fields should vanish
on the cathode.

The treatment of free boundary solutions (not necessarily being a
variational solution) leads to the notion of strong solutions with a free
boundary:

\begin{definition}
We say that $u\in C^{1}([0,1))\cap K$ is a strong solution with a free
boundary if there exists $\xi \in \lbrack 0,1)$ such that (\ref{flat
condition}) holds and, in addition, $\frac{1}{\sqrt{u}}\in L^{1}(\xi ,1)$,
and $u^{\prime \prime }(y)=\displaystyle\frac{j}{\sqrt{u(y)}}$ a.e. $y$ $\in
(\xi ,1).$
\end{definition}

Notice that, in this case, we must understand the differential equation as%
\begin{equation*}
-u^{\prime \prime }(y)+\displaystyle\frac{j}{\sqrt{u(y)}}\chi _{\left\{
u>0\right\} }=0\quad y\in (0,1),
\end{equation*}%
where $\chi _{\left\{ u>0\right\} }$ denotes the characteristic function of
the set $\left\{ y\in (0,1):u(y)>0\right\} $, since otherwise the singular
term is not well defined on $[0,\xi ]$.

\bigskip

The following result was communicated to me by Ha\"{\i}m Brezis on 2004. Although there are different similar treatments in
the literature (see, e.g., \cite{Degond-Raviart1991}) the uniqueness part
and the argument used in the proof is completely new.

\begin{theorem}
\label{thm Brezis}(H. Brezis) There exists a unique variational solution.
Moreover, if we define $j^{\ast }=\frac{4}{9},$ 

\noindent then:

\noindent a) if $j=j^{\ast }$ the variational solution is a flat solution
and it is given by 
\begin{equation}
u(y)=y^{4/3}.  \label{Solut *}
\end{equation}

\noindent b) if $j>j^{\ast }$the variational solution is a free boundary
solution with 
\begin{equation}
\xi =1-\frac{1}{\sqrt{\displaystyle\frac{9}{4}j}},
\end{equation}%
\noindent and it is given by $u(y)=A(y-\xi )_{+}^{4/3}$ \noindent with $A=(%
\frac{9}{4}j)^{2/3}.$ 

\noindent c) if \ $0<j<j^{\ast }$ then the variational solution is such that 
$u>0$ on $(0,1]$ and $u^{\prime }(0)=K_{0}>0,\quad $for some $K_{0}=K_{0}(j).
$ 
\end{theorem}

\begin{remark}\rm
In fact, from the proof we will see that if $j\rightarrow +\infty $ then $%
\xi \rightarrow 1$ (a sort of limit of \textquotedblleft boundary
layer\textquotedblright\ type). Moreover, if $j\rightarrow +0$ then $%
u(y)\rightarrow y$ in $C^{0}([0,1]).$
\end{remark}

\bigskip

\noindent \textit{Proof.} Let us start by proving the uniqueness of
variational free boundary type solutions. By well-known results, (\ref{e2})
admits always a minimizer which satisfies $u^{\prime \prime }=\frac{j}{\sqrt{%
u}}$ on the set $[\xi <y<1].$ Then, by multiplying by $u^{\prime }$ we get

\begin{equation}
\left\{ 
\begin{array}{l}
\displaystyle\frac{du}{dy}=\sqrt{4j}\left( u(y)\right) ^{1/4} \\ [.3cm]
u(1)=1,%
\end{array}%
\right.   \label{e8}
\end{equation}%
a problem with the uniqueness of solutions since $u^{1/4}$ is monotone
increasing. By Leibniz's formula for separable ordinary differential
equations, we get 
\begin{equation*}
u(y)=\left( \displaystyle\frac{3\sqrt{j}}{2}\right) ^{4/3}(y-\xi )_{+}^{4/3},
\end{equation*}
which proves a) and b). Let us consider now the case of $j\in (0,\frac{4}{9}%
).$ Note that if we arrive to prove that $\frac{du}{dy}(0)=K_{0}>0$ then
after multiplying by $u^{\prime }$ we will get that 
\begin{equation*}
u^{\prime }(y)=\displaystyle\sqrt{4j\sqrt{u(y)}+K_{0}^{2}},
\end{equation*}

\noindent and the Leibniz formula would lead to the implicit expression of
the solution 
\begin{equation}
\displaystyle\int_{0}^{u(y)}\frac{du}{\sqrt{4j\sqrt{u}+k_{0}^{2}}}=y.
\label{e14}
\end{equation}%
A more direct way to apply the above argument is the following: let

\begin{equation}
H(\sigma )=\displaystyle\int_{0}^{\sigma }\frac{ds}{\sqrt{s^{1/2}+1}},\quad 
\text{for }\sigma >0.
\end{equation}

\noindent Then $H(0)=0$ and

\begin{equation}
H^{\prime }(\sigma )=\displaystyle\frac{1}{\sqrt{{\sqrt{\sigma }+1}}}>0.
\end{equation}%
\noindent Thus $H$ can be inverted, $\sigma =H^{-1}(\tau )$, and we get that 
$H^{-1}(0)=0$ and

\begin{equation}
\displaystyle\frac{d}{d\tau }H^{-1}(\tau )=\sqrt{\sqrt{H^{-1}(\tau )}+1}>0.
\label{e18}
\end{equation}

\noindent Let us prove that any variational solution can be written in the
form $u(y)=\displaystyle\frac{1}{A}H^{-1}(By),$ for suitable constants $A$
and $B$. Note that%
\begin{equation}
u^{\prime }(y)=\displaystyle\frac{B}{A}\sqrt{\sqrt{Au(y)}+1}\quad \text{\
and thus }u^{\prime }(0)=\frac{B}{A}.  \label{e20}
\end{equation}

\noindent Indeed, since $H(Au(y))=By$ we get 
\begin{equation*}
H^{\prime }(Au)A\displaystyle\frac{du}{dy}=B,
\end{equation*}%
i.e., 
\begin{equation*}
\displaystyle\frac{1}{\sqrt{\sqrt{Au(y)}+1}}A\frac{du}{dy}=B,
\end{equation*}%
which implies \eqref{e20}. Moreover 
\begin{equation}
\begin{array}{ll}
\displaystyle\frac{d^{2}u}{dy^{2}}(y)&\displaystyle =\frac{B}{A}\frac{1}{2}\left( \sqrt{%
Au(y)}+1\right) ^{-1/2}\frac{1}{2}\left( Au(y)\right) ^{-1/2}A\displaystyle%
\frac{du}{dy}(y)\\ [.35cm]
& =\displaystyle\frac{B}{A}\displaystyle\frac{1}{2}\left( \sqrt{Au(y)}%
+1\right) ^{-1/2}\displaystyle\frac{1}{2}\left( Au(y)\right) ^{-1/2}B\sqrt{%
\sqrt{Au(y)}+1}  \\ [.35cm]
& =\displaystyle\frac{B^{2}}{4A}\displaystyle\frac{1}{\sqrt{Au(y)}}\\[.35cm]
& =\displaystyle\frac{B^{2}}{4A^{3/2}}\frac{1}{\sqrt{u(y)}}. 
\end{array}%
\end{equation}%
Thus, we must determine $A$ and $B$ by the conditions $\displaystyle\frac{%
B^{2}}{4A^{3/2}}=j$ and $u(1)=1.$ Then we get the condition 
\begin{equation}
H(A)=B=2\sqrt{j}A^{3/4},
\end{equation}%
which has a unique solution $A$ (and then a unique value of $B$) since $H$
is strictly increasing, $H(0)=0$ and $\displaystyle\lim_{\sigma \rightarrow
+\infty }H(\sigma )=+\infty $. Then 
\begin{equation*}
\displaystyle u^{\prime }(0)=\frac{B_{0}}{A_{0}}=\frac{H(A_{0})}{A_{0}},
\end{equation*}%
with $A_{0}>0$ solution of 
\begin{equation*}
\displaystyle\frac{H(A_{0})}{A_{0}^{3/4}}=2\sqrt{j}.
\end{equation*}%
The uniqueness of the constants $A$ and $B$ also proves the uniqueness of
the variational solution, and part c) follows.$\fin$

\begin{remark}\rm
\noindent \label{Remk Gradient estimate} In cases a) and
b) we can improve the regularity: we have that $u\in W^{2,p}(\xi ,1)\quad \forall p\in \lbrack 1,\frac{3}{2}),$ and in
case c) $u\in W^{2,p}(0,1)\quad \forall p\in \lbrack 1,2)$. In fact, we have
the \textit{sharp gradient estimate}%
\begin{equation}
\left\vert u^{\prime }(y)\right\vert \leq Cu^{1/4}(y)\text{, for some }C>0,%
\text{ for any }y\in (0,1).  \label{Estim gradient}
\end{equation}%
\noindent Estimates of this nature play an important role in the study of
the existence of solutions of the associated parabolic problem (the so-called
\textquotedblleft quenching problem\textquotedblright ): see, e.g., \cite%
{Phillips} and \cite{Dao-Diaz}, among others.
\end{remark}

\begin{remark}\rm
Note also that $j^{\ast }=\sup \{j>0:u^{\prime }(0)>0$ with $u$ variational
solution of (\ref{e2})$\}$. Something similar was already mentioned in the
space charge literature (see, e.g., \cite{Kirsten} p. 381 and \cite%
{Luginsland} p. 2372).
\end{remark}

\bigskip

\subsection{The distributed one-dimensional current density}

For different purposes, it is useful to consider a similar problem to the
above formulation but now for some $j=j(y)$ with $j\in L_{loc}^{1}(0,1),$ $%
j\geq 0$. The problem under consideration is 
\begin{equation}
\left\{ 
\begin{array}{l}
-u^{\prime \prime }(y)+\displaystyle\frac{j(y)}{\sqrt{u(y)}}\chi _{\left\{
u>0\right\} }=0\quad y\in (0,1),u\geq 0 \\ [.35cm]
u\geq 0\quad y\in (0,1), \\ [.2cm]
u(0)=0\quad u(1)=1,%
\end{array}%
\right.   \label{e35}
\end{equation}%
where $\chi _{\left\{ u>0\right\} }$ denotes again the characteristic
function of the set $\left\{ y\in (0,1):u(y)>0\right\} $. Since the
non-linear term is decreasing, we will apply the method of super and
subsolutions.

\begin{definition}
\label{def9} Given $p\geq 1,$ a function $\overline{u}\in H^{1}(0,1)$ is a $%
p-$strong supersolution of problem \eqref{e35}, if 

i) $\overline{u}\geq 0$ and $\displaystyle\frac{j}{\sqrt{\overline{u}}}\in
L^{p}(0,1),$ 

ii) $\overline{u}^{\prime \prime }\in L^{1}(0,1),$

iii) $-\overline{u}^{\prime \prime }+\displaystyle\frac{j}{\sqrt{\overline{u}%
}}\chi _{\left\{ \overline{u}>0\right\} }\geq 0$ a.e. in $(0,1),$

iv) $\overline{u}(0)=0$ and $\overline{u}(1)\geq 1.$
\end{definition}

The notion of $p-$strong subsolution \underline{$u$}$\in H^{1}(0,1)$ is
similar, by replacing the inequalities $\geq $ \ in iii) and iv) by the
inequality $\leq $, but always with $\underline{u}(0)=0$. We recall here a variation of the iterative method of
super and subsolutions (see, e.g., the exposition and general ideas in the
monograph \cite{P1992}) which applies to our framework with a singular
absorption. The proof is a special case of a more general statement, which
will be given later (see Theorem \ref{Thm Existence super-subs} and Remark \ref{super subs very weak}).

\begin{theorem}
\label{t10} Let $j\in L_{loc}^{1}(0,1),$ $j\geq 0.$ Assume that there exists 
$p>1,$ a $p-$strong supersolution $u^{0}$ and a $p-$strong subsolution $u_{0}
$, such that 
\begin{equation}
0<u_{0}(y)\leq u^{0}(y)\text{ a.e. }y\in (0,1).  \label{e36}
\end{equation}%
Then, problem \eqref{e35} has a maximal solution $u^{\ast }$ and a minimal
solution $u_{\ast }$ on the \textquotedblleft interval\textquotedblright 
\hspace{0mm} $[u_{0},u^{0}]$ of $H^{1}(0,1)$, i.e., such that 
\begin{equation}
u_{0}\leq u_{\ast }\leq u^{\ast }\leq u^{0}\text{ a.e. }y\in (0,1).
\label{e37}
\end{equation}%
Moreover $u_{\ast }^{\prime \prime }\text{, }u^{\prime \prime \ast }\in
L^{p}(0,1)$.$\fin$
\end{theorem}

\bigskip

We can consider some special cases of $j(y)$ that allow to get some results
in the line of Theorem \ref{thm Brezis}: 
\begin{equation}
j(y)=\lambda y^{q}\text{ with }q\in (-\frac{1}{2},1),\text{ }\lambda >0.
\label{Assump j(x)}
\end{equation}%
Note that, obviously, $q=0$ corresponds to the case treated before and that
our study will consider cases in which $j(0)=0$ ($q\in (0,1)$) as well as
cases in which $j(0)=+\infty $ ($q\in (-\frac{1}{2},0)$).

\bigskip

\begin{theorem}
\label{TheoremD-2022}Assume (\ref{Assump j(x)}). Then if we define $\lambda
_{q}^{\ast }=\frac{2(1+2q)(2+q)}{9}$ we have:

\noindent a) if $\lambda =$ $\lambda _{q}^{\ast }$ the function $%
u(y)=y^{(4+2q)/3}$ is a flat solution of the problem.

\noindent b) if $\lambda >$ $\lambda _{q}^{\ast }$ the function $%
u(y)=A(y-\xi )_{+}^{(4+2q)/3}$ is a free boundary solution with $\displaystyle \xi =1-%
\frac{1}{\left[ \frac{\lambda }{\lambda _{q}^{\ast }}\right] ^{1/(2+q)}},$
and $A=(\frac{\lambda }{\lambda _{q}^{\ast }})^{2/3}.$ Moreover,%
\begin{equation*}
\frac{j(y)}{\sqrt{u(y)}}\chi _{\left\{ u>0\right\} }=\frac{\lambda }{\sqrt{A}%
}(y-\xi )^{-\frac{(2-2q)}{3}}\chi _{\left\{ y>\xi \right\} }\in L^{1}(0,1).
\end{equation*}

\noindent c) if $\ 0<\lambda <$ $\lambda _{q}^{\ast },$ then the unique
solution $u$ is such that $u>0$ on $(0,1]$ and $u^{\prime }(0)=K_{0}>0,\quad 
$for some $K_{0}=K_{0}(\lambda ).$ 
\end{theorem}

\begin{remark}\rm
\noindent\ If we use a different notation, $\beta =-q$ and 
$\alpha =(4-2\beta )/3$, then we get the existence of solutions of
the form $u(y)=y^{\alpha }$, with $\alpha \in (1,4/3)$,
once we assume $\beta \in (0,1/2)$. Curiously enough, the
constraint $3\alpha /2+\beta =2$ arises also in Theorem \ref{TheomMain} dealing with the two-dimensional problem (see Section 4). 
\end{remark}

The key idea of the proof is to build a family of explicit solutions of the
differential equation inspired by Chapter 2 of \cite{D}.

\begin{lemma}
\label{l11} Let $u_{m}(y)=Cy^{\frac{2}{1+m}}$, with $C>0$ and $m\in (0,1)$.
\noindent Then 
\begin{equation}
-u_{m}^{\prime \prime }(y)+\displaystyle\frac{j_{m}(y)}{\sqrt{u_{m}(y)}}=0,
\label{e48}
\end{equation}%
where 
\begin{equation}
j_{m}(y)=\displaystyle\frac{2C\sqrt{C}(1-m)}{(1+m)^{2}}y^{\frac{(1-2m)}{(1+m)%
}}\quad y\in (0,1).  \label{e49}
\end{equation}
\end{lemma}

\noindent \textit{Proof.} It suffices to check that 
\begin{equation*}
(u_{m})^{\prime }(y)=C\displaystyle\frac{2}{1+m}y^{\left( \frac{2}{1+m}%
-1\right) }=\frac{2C}{1+m}y^{\frac{1-m}{1+m}},
\end{equation*}%
\begin{equation}
(u_{m})^{\prime \prime }(y)=\displaystyle\frac{2C(1-m)}{(1+m)^{2}}y^{\left( 
\frac{1-m}{1+m}-1\right) }=\frac{2C(1-m)}{(1+m)^{2}}y^{\frac{-2m}{1+m}}.
\label{e50}
\end{equation}%
Since $\sqrt{u_{m}(y)}=\sqrt{C}y^{\frac{1}{1+m}},$ we get %
\begin{equation*}
-(u_{m})^{\prime \prime }(y)=\displaystyle\frac{-j_{m}(y)}{\sqrt{u_{m}(y)}}=%
\frac{-1}{\sqrt{C}}y^{\frac{-1}{1+m}}j_{m}(y).
\end{equation*}%
Applying \eqref{e50} we arrive at (\ref{e48}) with $j_{m}(y)$ given by %
\eqref{e49}.$\fin$

\bigskip

\begin{corollary}
\label{Corol explicit j(x)}Given $j(y)=\lambda y^{q}$ with $q\in (-\frac{1}{2%
},1)$ and $\lambda >0,$ the function $u_{q}(y)=C_{q}(\lambda )y^{\frac{4+2q}{%
3}},$with 
\begin{equation}
C_{q}(\lambda )=\left[ \frac{9\lambda }{2(1+2q)(2+q)}\right] ^{2/3},
\label{Constant Change exponents}
\end{equation}%
\noindent is a solution of 
\begin{equation*}
-u_{q}^{\prime \prime }+\displaystyle\frac{\lambda y^{q}}{\sqrt{u_{q}}}=0,%
\text{ in }(0,1).
\end{equation*}
\end{corollary}

\noindent \textit{Proof.} We impose $\frac{(1-2m)}{(1+m)}=q$ \ and 
\begin{equation*}
\frac{2C\sqrt{C}(1-m)}{(1+m)^{2}}=\lambda .
\end{equation*}%
Then we get $m=\frac{1-q}{2+q}$ and $C(\lambda )$ given by (\ref{Constant
Change exponents}). Then the conclusion is a direct consequence of Lemma~\ref{l11}.$\fin$

\bigskip

\noindent \textit{Proof of Theorem \ref{TheoremD-2022}. }Conclusions a) and
b) are a direct consequence of Corollary \ref{Corol explicit j(x)}. To prove
part c) it suffices to build a subsolution $\underline{u}(y)\leq y$
satisfying the searched properties, since $\overline{u}(y)=y$ is a
supersolution satisfying that $\underline{u}(y)\leq \overline{u}(y)$ a.e. $%
y\in (0,1),$ and then, by Theorem \ref{t10}, there exists a minimal solution $%
u_{\ast }$ of (\ref{e35}) satisfying 
\begin{equation*}
\underline{u}(y)\leq u_{\ast }(y)\leq \overline{u}(y)\text{ a.e. }y\in (0,1).
\end{equation*}%
Then if $\underline{u}^{\prime }(0)>0$ we conclude, as desired, that $%
u_{\ast }^{\prime }(0)>0.$

\noindent A construction of a suitable subsolution is the following: let $%
\varepsilon >0$ be small enough, to be determined later, let $%
u_{\#}(y)=y^{(4+2q)/3}$ and define %
\begin{equation*}
\underline{u}(y)=\frac{u_{\#}(y)+\varepsilon y}{1+\varepsilon }.
\end{equation*}%
Then it is clear that $\underline{u}(0)=0$, $\underline{u}(1)=1$ and since $%
\underline{u}^{\prime }(y)=\frac{u_{\#}^{\prime }(y)+\varepsilon }{%
1+\varepsilon }$ we get $\underline{u}^{\prime }(0)=\frac{\varepsilon }{%
1+\varepsilon }>0.$ Moreover, 
\begin{equation*}
\underline{u}^{\prime \prime }(y)=\frac{u_{\#}^{\prime \prime }(y)}{%
1+\varepsilon }=\frac{\lambda _{q}^{\ast }y^{q}}{(1+\varepsilon )\sqrt{%
u_{\#}(y)}}.
\end{equation*}%
Then, concerning the differential equation, we can have the condition for
subsolution%
\begin{equation*}
-\underline{u}^{\prime \prime }(y)+\displaystyle\frac{\lambda y^{q}}{\sqrt{%
\underline{u}(y)}}=-\frac{\lambda _{q}^{\ast }y^{q}}{(1+\varepsilon )\sqrt{%
u_{\#}(y)}}+\displaystyle\frac{\sqrt{1+\varepsilon }\lambda y^{q}}{\sqrt{%
u_{\#}(y)+\varepsilon y}}\leq 0
\end{equation*}%
if we have%
\begin{equation}
\frac{\sqrt{1+\varepsilon }\lambda }{\sqrt{u_{\#}(y)+\varepsilon y}}\leq 
\frac{\lambda _{q}^{\ast }}{(1+\varepsilon )\sqrt{u_{\#}(y)}}.
\label{searched condition}
\end{equation}%
But, obviously
\begin{equation*}
\frac{\lambda _{q}^{\ast }}{(1+\varepsilon )\sqrt{u_{\#}(y)+\varepsilon y}}%
\leq \frac{\lambda _{q}^{\ast }}{(1+\varepsilon )\sqrt{u_{\#}(y)}}.
\end{equation*}%
Then we arrive at the searched inequality (\ref{searched condition}) if we
have%
\begin{equation*}
\lambda \leq \frac{\lambda _{q}^{\ast }}{(1+\varepsilon )^{3/2}}.
\end{equation*}%
This is clearly true since $\lambda <\lambda _{q}^{\ast }$: it suffices to
take $\varepsilon >0$ such that 
\begin{equation*}
\varepsilon <\left( \frac{\lambda _{q}^{\ast }}{\lambda }\right)
^{2/3}-1.
\end{equation*}
\fineq

\bigskip

The functions $u_{q}(y)$ given in Corollary \ref{Corol explicit j(x)} can be
used as super and subsolutions to get the existence of a flat solution for a
more general function $j(y)$. The result that follows is just a sample of
many other possible results.

\begin{theorem}
\label{t115} Let $-\frac{1}{2}<r<q<0.$ Assume 
\begin{equation*}
\lambda _{r}y^{r}\geq j(y)\geq \lambda _{q}y^{q}\text{ }\quad \text{a.e. }%
y\in (0,1),
\end{equation*}%
with 
\begin{equation}
\lambda _{q}=\lambda _{q}^{\ast }\text{ and }\lambda _{r}>\lambda _{r}^{\ast
}\text{ with }\lambda _{r}\text{ large enough,}
\label{possible freeboundary}
\end{equation}%
or%
\begin{equation}
\lambda _{q}<\lambda _{q}^{\ast }\text{ and }\lambda _{r}=\lambda _{r}^{\ast
}\text{ with }\lambda _{q}\geq 0\text{ small enough.}
\label{case positive subsol}
\end{equation}

\noindent Then there exists a maximal solution $u^{\ast }$ and a minimal
solution  $u_{\ast }$ of problem (\ref{e35}). Moreover, in case of (\ref%
{possible freeboundary}) 
\begin{equation*}
A_{r}(y-\xi _{r})_{+}^{\frac{4+2r}{3}}\leq u_{\ast }(y)\leq u^{\ast }(y)\leq
y^{\frac{4+2q}{3}}\text{ a.e. }y\in (0,1),
\end{equation*}%
for some $A_{r}>0$ and $\xi _{r}\in (0,1)$, and in case (\ref{case positive
subsol})%
\begin{equation*}
y^{\frac{4+2r}{3}}\leq u_{\ast }(y)\leq u^{\ast }(y)\leq u_{q}(y)\text{ a.e. 
}y\in (0,1),
\end{equation*}%
with $u_{q}(y)$ the solution of problem (\ref{e35}) corresponding to $%
j(y)=\lambda _{q}y^{q}$ (when $\lambda _{q}\in (0,\lambda _{q}^{\ast })$)
mentioned in part c) of Theorem \ref{TheoremD-2022}. Moreover if $\lambda
_{q}=0$ we can take as $u_{q}$ the function $u_{q}(y)=y.$
\end{theorem}

\noindent \textit{Proof.} We recall that for $a>b\geq 0$ we have $%
0<y^{a}<y^{b}<1$ for any $y\in (0,1).$ Thus, equivalently, if $-a<-b\leq 0$
then $y^{-a}>y^{-b}>1$ for any $y\in (0,1).$

\noindent We take now $-a=r$ and $-b=q$ so $y^{r}>y^{q}$. The conclusion
will be obtained through the application of Theorem \ref{TheoremD-2022} to
different choices of the super and subsolution. In case (\ref{possible
freeboundary}) we take as supersolution $u^{0}(y)=y^{\frac{4+2q}{3}}$ since 
\begin{equation*}
-u^{0\prime \prime }(y)+\displaystyle\frac{j(y)}{\sqrt{u^{0}(y)}}\geq
-u^{0\prime \prime }(y)+\displaystyle\frac{\lambda _{q}^{\ast }y^{q}}{\sqrt{%
u^{0}(y)}}=0.
\end{equation*}%
On the other hand, as a subsolution we can take $u_{0}(y)=A_{r}(y-\xi
_{r})_{+}^{\frac{4+2r}{3}}$ for some $A_{r}>0$ and $\xi _{r}\in (0,1)$ as
indicated in part b) of Theorem \ref{TheoremD-2022} (remember that $\lambda
_{r}>\lambda _{q}^{\ast }$) since 
\begin{equation*}
-u_{0}^{\prime \prime }(y)+\displaystyle\frac{j(y)}{\sqrt{u_{0}(y)}}\leq
-u_{0}^{\prime \prime }(y)+\displaystyle\frac{\lambda _{r}y^{r}}{\sqrt{%
u_{0}(y)}}=0.
\end{equation*}

\noindent We have $u_{0}(0)=u^{0}(0)=0,$ $u_{0}(1)=u^{0}(1)=1$. Moreover, if 
$\lambda _{r}$ is large enough we have 
\begin{equation*}
u_{0}(y)=A_{r}(y-\xi _{r})_{+}^{\frac{4+2r}{3}}<y^{\frac{4+2q}{3}}=u^{0}(y)%
\text{ for any }y\in (0,1),
\end{equation*}%
(in spite that $0<\frac{4+2r}{3}<\frac{4+2q}{3}$ and thus $y^{\frac{4+2q}{3}%
}<y^{\frac{4+2r}{3}}$ for $y\in (0,1)$). Indeed, we can take $A_{r}>0$ and $%
\xi _{r}\in (0,1)$ so that $u_{0}^{\prime }(1)>u^{0^{\prime }}(1)=\frac{4+2q%
}{3}$ and thus both functions are well-ordered.

\noindent In case (\ref{case positive subsol}) we take as subsolution the
function $u_{0}(y)=y^{\frac{4+2r}{3}}$ and as supersolution the function $%
u^{0}(y)=u_{q}(y)$ the solution of problem (\ref{e35}) corresponding to $%
j(y)=\lambda _{q}y^{q}$ when $\lambda _{q}\in (0,\lambda _{q}^{\ast })$
mentioned in part c) of Theorem \ref{TheoremD-2022}. The conditions in terms
of the differential equation are satisfied, as before. Moreover, again, $%
u_{0}(0)=u^{0}(0)=0,$ $u_{0}(1)=u^{0}(1)=1$, and the inequality $u_{0}(y)=y^{%
\frac{4+2r}{3}}<u_{q}(y)=u^{0}(y)$ for any $y\in (0,1),$ holds once that $%
\lambda _{q}>0$ is small enough (since we know that $u_{q}(y)\rightarrow y$
when $\lambda _{q}\rightarrow 0$). If $\lambda _{q}=0$ the function $%
u_{q}(y)=y$ is a supersolution and the proof ends.$\fin$

\begin{remark}\rm
\noindent The gradient estimate mentioned when $j$ is a
constant is no longer valid when $j(y)=\lambda y^{-\beta }$ with $%
\beta \in (0,\frac{1}{2})$. If $\lambda \geq \lambda _{\beta }$, for some $\lambda _{\beta }>0$
\begin{equation}
\left\vert u^{\prime }(y)\right\vert \leq Cu^{\frac{1+2\beta }{4+2\beta }}(y),
\quad \text{for some $C>0$, for any $y\in (0,1)$}.
\end{equation}
\end{remark}

\section{The 2-d parallel-plate geometry}

As mentioned in the Introduction, Theorem \ref{TheomMain} makes precise, a
conjecture by A. Rokhlenko \cite{Rokhlenko}: if $j(x)$ behaves as $%
A/\left\vert x\right\vert ^{\beta },$ for some $\beta \in (0,1/2)$, for $%
x\in (-a,0)$ then there exists a weak solution $u(x,y)$ of (\ref{BVP}), and $%
u$ \textquotedblleft behaves\textquotedblright\ (near the cathode $%
[-a,a]\times \{0\}$) as $y^{\alpha }$ for some $\alpha \in (1,4/3),$ with $%
\alpha =4/3-2\beta $. \ We will obtain other auxiliary results for the more
general case of $j:(-a,b)\mathbb{\rightarrow \lbrack }0,+\infty )$ such that 
\begin{equation}
\left\{ 
\begin{array}{lr}
j(x)>0 & \text{if }x\in (-a,0),\text{ }j\in L_{loc}^{1}(-a,0), \\ 
j(x)=0 & \text{if }x\in (0,b),%
\end{array}%
\right.  \label{Hypo j a b}
\end{equation}
with a possible singularity of the current density function $j(x)$ only
expected at the origin $x=0$. Our main goal is, therefore, to get the
transition, near $x=0$, between flat and linear profiles $u(.,y).$

As mentioned before, since in other frameworks the solutions may vanish in
some parts of the spatial domain (\textit{dead cores}),  we can reformulate
the PDE as%
\begin{equation}
P_{a,b,j}=\left\{ 
\begin{array}{lr}
-\Delta u+\frac{j(x)}{\sqrt{u}}\chi _{\{u>0\}}=0 & x\in (-a,b),\text{ }y\in
(0,1), \\ [.3cm]
u(x,0)=0 & x\in (-a,b), \\ [.15cm]
u(x,1)=1 & x\in (-a,b), \\ [.15cm]
u(-a,y)=y^{4/3} & y\in (0,1), \\ [.15cm]
u(b,y)=y^{{}} & y\in (0,1).%
\end{array}%
\right. 
\end{equation}

This Section will be structured in several subsections. In subsection 4.1 we
present the adaptation of the super and subsolutions method to the problem $%
P_{a,b,j}$. The construction of \textit{positive subsolutions }satisfying
the additional conditions $AC_{a,b}$ will be presented in subsection 4.2.
The proof is made by constructing suitable super and subsolutions for
several auxiliary problems and matching suitably the corresponding solutions. In particular, we will use the global bifurcation diagram (in terms of
the parameter $\lambda $) associated to the auxiliary problem (\ref{Problem
Psi}). In subsection 4.3 we will end the proof of Theorem \ref{TheomMain} by
constructing a suitable supersolution which is flat on $(-a,0)\times \{0\}$.
Finally, in subsection 4.4 we will prove the uniqueness of non-degenerate
solutions.

\subsection{On the existence of solutions of $P_{a,b,j}$ via super and
subsolutions method}

Due to the singularity of the nonlinear term, the super and subsolutions
method needs to be applied under some adequate conditions. Let $\Omega
=(-a,b)\times (0,1)$ and $j$ as in (\ref{Hypo j a b}).

\begin{definition}
\textrm{\textit{A function }$u^{0}\in W^{1,1}(\Omega )$, is said a $p-$%
positive supersolution of $P_{a,b,j}$ if $\ u^{0}\geq 0$, $\frac{j}{\sqrt{%
u^{0}}}\in L^{p}(\Omega )$, for some $p\geq 1,$ and it verifies in a very
weak sense that 
\begin{equation*}
\left\{ 
\begin{array}{lr}
-\Delta u^{0}\geq -\frac{j(x)}{\sqrt{u^{0}}} & \text{in }\Omega , \\[0.3cm]
u^{0}(x,0)\geq 0 & x\in (-a,b), \\[0.15cm]
u^{0}(x,1)\geq 1 & x\in (-a,b), \\[0.15cm]
u^{0}(-a,y)\geq y^{4/3} & y\in (0,1), \\[0.15cm]
u^{0}(b,y)\geq y & y\in (0,1).%
\end{array}%
\right. 
\end{equation*}%
The notion of $p-$positive subsolution $u_{0}$ is introduced similarly,
i.e.,  $u_{0}\in W^{1,1}(\Omega )$ with $u_{0}> 0$ and $\frac{j}{\sqrt{%
u_{0}}}\in L^{p}(\Omega )$, satisfies that
\begin{equation*}
\left\{ 
\begin{array}{lr}
-\Delta \mathrm{u_{0}}\leq -\frac{j(x)}{\sqrt{\mathrm{u_{0}}}} & \text{in }%
\Omega , \\[0.3cm]
\mathrm{u_{0}}(x,0)=0 & x\in (-a,b), \\[0.15cm]
\mathrm{u_{0}}(x,1)\leq 1 & x\in (-a,b), \\[0.15cm]
\mathrm{u_{0}}(-a,y)\leq y^{4/3} & y\in (0,1), \\[0.15cm]
u^{0}(b,y)\leq y & y\in (0,1).%
\end{array}%
\right. 
\end{equation*}%
}
\end{definition}

\begin{theorem}
\label{Thm Existence super-subs}\textit{Assume that }%
\begin{equation}
\left\{ 
\begin{array}{c}
\text{\textit{there\ exists\ }}p>1\text{, \textit{a\ }}\mathit{p-}\text{%
\textit{positive supersolution\ }}\ u^{0}\text{ \textit{and a }}\mathit{p-}%
\text{positive \textit{subsolution }}\mathit{u}_{0} \\ 
\text{of }P_{a,b,j}\text{\textit{\ such\ that\ }}u_{0}\leq u^{0}\text{%
\textit{\ a.e. in\ }}\Omega .%
\end{array}%
\right.   \label{Super_sub}
\end{equation}%
\textit{Then problem }$P_{a,b,j}$\textit{\ possesses a minimal and a maximal
solutions }$u_{\ast }$\textit{\ and }$u^{\ast }$\textit{\ in the interval }$%
[u_{0},u^{0}]$\textit{, i.e.,} 
\begin{equation*}
u_{0}\leq u_{\ast }\leq u^{\ast }\leq u^{0}\text{ a.e.in }\Omega .
\end{equation*}
\end{theorem}

\noindent \textit{Proof.}\textbf{\ }We define the iterative schemes
(starting with $u^{0}$ and $u_{0}$) 
\begin{equation*}
\left\{ 
\begin{array}{lr}
-\Delta u^{n}=-\frac{j(x)}{\sqrt{u^{n-1}}} & \text{in }\Omega , \\ [.3cm]
u^{n}(x,0)=0 & x\in (-a,b), \\ [.15cm]
u^{n}(x,1)=1 & x\in (-a,b), \\ [.15cm]
u^{n}(-a,y)=y^{4/3} & y\in (0,1), \\ [.15cm]
u^{n}(b,y)=y & y\in (0,1),%
\end{array}%
\right.
\end{equation*}%
and

\begin{equation*}
\left\{ 
\begin{array}{lr}
-\Delta u_{n}=-\frac{j(x)}{\sqrt{u_{n-1}}} & \text{in }\Omega , \\ [.3cm]
u_{n}(x,0)=0 & x\in (-a,b), \\ [.15cm]
u_{n}(x,1)=1 & x\in (-a,b), \\ [.15cm]
u_{n}(-a,y)=y^{4/3} & y\in (0,1), \\ [.15cm]
u_{n}(b,y)=y & y\in (0,1).%
\end{array}%
\right. 
\end{equation*}%
By using the comparison principle for the Laplace operator, we get that 
\begin{equation*}
0<u_{0}\leq u_{1}\leq ...\leq u_{n}\leq ...\leq u^{n}\leq ...\leq u^{1}\leq
u^{0}\text{ \ a.e. in \ }\Omega ,
\end{equation*}%
and so the sequences $\{u^{n}\},\{u_{n}\}$ converge (monotonically) in $%
L^{p}(\Omega )$ to some functions $u_{\ast }$\textit{\ and }$u^{\ast }$ and
the sequences $\{\frac{j}{\sqrt{u^{n-1}}}\}$, $\{\frac{j}{\sqrt{u_{n-1}}}\}$
are bounded in $L^{p}(\Omega )$ and converge also (monotonically) in $%
L^{p}(\Omega )$.$_{\blacksquare }$

\bigskip

\begin{remark}\rm 
\noindent \label{super subs very weak} With some minor modifications, the above result holds when the super and subsolutions are in the weighted space $\frac{j}{\sqrt{u^{0}}}\in
L^{p}(\Omega ,\delta )$, for some $p\geq 1,$ with $\delta =d((x,y),\partial
\Omega ).$ See, e.g., \cite{Brezis-Cazenave-Martel}, \cite{Souplet}, \cite%
{Diaz-Rakotoson} and \cite{montenegro2008sub} among many other papers (see also some applications in \cite{Brezis-Diaz-eds}). This allows a greater generality to treat the
singular term: for instance, the function $u^{0}(x,y)=y$ is such that $\frac{%
1}{\sqrt{u^{0}}}\in L^{p}(\Omega )$ for $p\in \lbrack 1,2),$ $\frac{1}{\sqrt{%
u^{0}}}\notin L^{2}(\Omega )$ but $\frac{1}{\sqrt{u^{0}}}\in L^{2}(\Omega
,\delta ).$
\end{remark}
\begin{remark}\rm
\noindent \label{Rmk Matching picos}By some standard approximating
arguments, it is well-known that the above notion of $p-$super and
subsolutions of $P_{a,b,j}$ can be extended to the case in which the
diffusion term generates an additional distribution over a simple curve $%
\Gamma $ separating $\Omega $ in two different parts (matching without a $%
W^{1,p}$-contact: see \cite{Ilin-Kalash-Oleinik} and  \cite{Berest-Lions}). Then, for instance, the
subsolution $u_{0}$ of problem $P_{a,b,j}$ is allowed to satisfy 
\begin{equation*}
-\Delta u_{0}+\frac{j(x)}{\sqrt{u_{0}}}\leq 0 \text{ in }%
\mathcal{D}^{\prime }(\Omega ).
\end{equation*}
\end{remark}

The existence of a (not-flat) supersolution $u^{0}$ can be easily proved.

\begin{lemma}
\textbf{\label{Lemm supersolution y}} \textit{Let }$j(x)$ \textit{satisfying
(\ref{Hypo j a b}). }Then\textit{\ the function }$u^{0}(x,y)=y$\textit{\ is
a positive }$p-$\textit{supersolution of }$P_{a,b,j}$ for any $p\in \lbrack
1,2)$\textit{.}
\end{lemma}

\noindent \emph{Proof.} It is a trivial fact since 
\begin{equation*}
-\Delta u^{0}=0\geq -\frac{j(x)}{\sqrt{u^{0}}}\text{ \ }x\in (-a,b),\text{ }%
y\in (0,1),
\end{equation*}%
and%
\begin{equation*}
\left\{ 
\begin{array}{lr}
u^{0}(x,0)=0 & x\in (-a,b), \\ [.15cm]
u^{0}(x,1)=1 & x\in (-a,b), \\ [.15cm]
u^{0}(-a,y)\geq y^{4/3} & y\in (0,1), \\ [.15cm]
u^{0}(b,y)=y & y\in (0,1).%
\end{array}%
\right. 
\end{equation*}%
Moreover, $\frac{j(x)}{\sqrt{u^{0}}}\in L^{p}(\Omega )$ for any $p\in
\lbrack 1,2)$.$\fin$

\bigskip

The construction of a positive subsolution is a very delicate task, which
will be presented in the following subsection. Before proving Theorem \ref%
{TheomMain} we can get an existence and uniqueness result of a positive
solution under a more general assumption on $j(x)$ than in Theorem \ref%
{TheomMain} but without ensuring that it is flat on $(-a,0)\times \{0\}.$ We
define

\begin{equation*}
\delta (x,y):=\mathrm{dist}((x,y),\partial \Omega ),
\end{equation*}%
which sometimes we shall denote simply as $\delta $.

\begin{theorem}
\label{Theorem No main} \textit{Assume that }$j(x)$\textit{\ satisfies }%
\begin{equation*}
j(x)=\frac{A}{(-x)^{\beta }},\text{ \textit{for} }x\in (-a,0)\text{
and }j(x)=0\text{ i\textit{f} }x\in (0,b),
\end{equation*}%
\textit{with}%
\begin{equation}
0<\beta <1/2\text{, }\beta \text{ and }A>0\text{ small\ enough}\mathit{.}
\label{Hipo beta}
\end{equation}%
\textit{\ Then there exists a weak solution} $u\in L^{2}(\Omega ;\delta )$ 
\textit{of} 
\begin{equation}
P_{a,b,j}=\left \{ 
\begin{array}{cc}
-\Delta u+\frac{j(x)}{\sqrt{u}}=0 & x\in (-a,b),\text{ }y\in (0,1), \\[0.3cm]
u(x,0)=0 & x\in (-a,b), \\[0.15cm]
u(x,1)=1 & x\in (-a,b), \\[0.15cm]
u(-a,y)=y^{4/3} & y\in (0,1), \\[0.15cm]
u(b,y)=y^{{}} & y\in (0,1),%
\end{array}%
\right. 
\end{equation}%
\textit{such that} 
\begin{equation}
C\delta (x,y)^{4/3}\leq u(x,y)\leq y\text{ a.e. }(x,y)\in (-a,b)\times (0,1),
\label{Estimate nondegenerate}
\end{equation}%
\textit{for some} $C>0$, and $Cy^{\alpha }\leq u(0,y),$ $y\in (0,1),$ with $%
\alpha =\frac{2}{3}(2-\beta )$.
\fineq
\end{theorem}

\bigskip

The solutions \ satisfying inequalities of the type (\ref{Estimate
nondegenerate}) are called \textit{nondegenerate} solutions. The
uniqueness of solutions given in Theorem \ref{Theorem No main} is a consequence of the
techniques introduced in the paper \cite{DiGiacomoni}. Indeed, let $\nu \in
(0,4/3],$ and define the class of functions%
\begin{equation}
\mathcal{M}(\nu ):=\Big\lbrace u\in L^{2}(\Omega ;\delta )\;\big\vert\;\text{%
such\ that }u(x,y)\geq C\delta (x,y)^{\nu }\quad \text{in }\Omega \text{,
for some }C>0\Big\rbrace.
\end{equation}%
We have

\begin{theorem}
\label{Thm uniqueness} \textit{Assume} $j(x)$ \textit{as in Theorem \ref%
{Theorem No main}. Then, there exists at most a solution} $u\in \mathcal{%
M}(\nu )$ \textit{of} $P_{a,b,j}.\fin$
\end{theorem}

\bigskip

The proof will be obtained through some smoothing estimates for some suitable
parabolic problem (see subsection 4.4). Finally, the proof of the
main result of this paper (Theorem \ref{TheomMain}) will be a direct
consequence of the method of super and subsolutions, and the following result
on partially flat supersolutions:

\begin{theorem}
\label{Theorem Supersolution}\textit{There exist} $A_{0},b_{0}>0$\textit{\
and }$\beta _{0}\in (0,\frac{1}{2}$\textit{) such that, if }$b\geq b_{0}>0,$ 
\textit{and if we assume} (\ref{Hypo j supersolution}) \textit{then there
exists a partially flat supersolution }$\overline{u}(x,y)$\textit{\ of
problem }$P_{a,b,j},$\textit{\ i.e., such that} $\overline{u}\in
L^{2}(\Omega ;\delta )$, 
\begin{equation}
\left \{ 
\begin{array}{cc}
-\Delta \overline{u}+\frac{j(x)}{\sqrt{\overline{u}}}\geq 0 & x\in (-a,b),%
\text{ }y\in (0,1), \\[0.3cm]
\overline{u}(x,0)\geq 0 & x\in (-a,b), \\[0.15cm]
\overline{u}(x,1)\geq 1 & x\in (-a,b), \\[0.15cm]
\overline{u}(-a,y)\geq y^{4/3} & y\in (0,1), \\[0.15cm]
\overline{u}(b,y)\geq y^{{}} & y\in (0,1),%
\end{array}%
\right. 
\end{equation}%
\textit{and } 
\begin{equation}
u(x,y)\leq \overline{u}(x,y)\leq C\delta (x,y)^{\overline{\alpha 
}}\text{ a.e. }(x,y)\in (-a,0)\times (0,1),  \label{Estimate supersolution}
\end{equation}%
\textit{for some} $C>0$, \textit{with} $\overline{\alpha }\in (1,\frac{4}{3})
$ and $u(x,y)$ the solution given in Theorem \ref{Theorem No
main}. \fineq
\end{theorem}

\begin{remark}\rm 
In order to have the zero flux condition on $(-a,0)$ we need some special
unbounded behaviour on $j(x)$ near $x=0.$ Given $x_{0}\in (-a,0)$ we can
construct (as in \cite{BBF}, \cite{D-Hern}, \cite{D}) a \textquotedblleft
barrier function\textquotedblright\ of the form%
\begin{equation}
\overline{u}(x,y:x_{0})=K\left\Vert (x,y)-(x_{0},0)\right\Vert
^{4/3}=K[(x-x_{0})^{2}+y^{2}]^{2/3},  \label{U-barre}
\end{equation}%
for a suitable $K>0$ in such a way that $\overline{u}$ let a local
supersolution once we assume%
\begin{equation*}
j(x)\geq C(-x)_{+}^{-\delta },\text{ }x\in (-a,a)
\end{equation*}%
for some suitable $\delta >0.$More precisely we can assume%
\begin{equation}
j(x)\geq 4\max (1,\frac{\sqrt{2}}{\sqrt{-x}}),\text{ }x\in (-a,0)\text{ \
and }j(x)=0\text{ for }x\in (0,a),  \label{singular j}
\end{equation}%
and that 
\begin{equation*}
j(x)[(x-x_{0})^{2}+y^{2}]^{-1/3}\in L^{p}(\Omega _{x_{0}}),\text{ for some }%
p>1,
\end{equation*}%
where $\Omega _{x_{0}}=\{(x,y)\in (-a,0)\times
(0,1):(x-x_{0})^{2}+y^{2}<(-x_{0})^{2}\}$, for any $x_{0}\in
(-a,0)$. Assume also that there exists  a p-nonnegative subsolution ${u}_{0}$ such that  
$u_{0}(x,y)\leq y$ a.e. in $\Omega$. Then, if $u$ satisfies the problem in the interval
$[u_{0},y]$, i.e., such that
\begin{equation*}
u_{0}\leq u\leq y\text{ \ \ }a.e.in\text{ }\Omega 
\end{equation*}%
we have that
\begin{equation}
u(x,y)\leq \max (1,\frac{2^{1/3}}{(-x_{0})^{1/3}}%
)[(x-x_{0})^{2}+y^{2}]^{2/3},  \label{estimate}
\end{equation}%
for any $(x,y)\in (-a,0)\times (0,1)$ and $x_{0}\in (-a,0)
$ such that $(x-x_{0})^{2}+y^{2}<(-x_{0})^{2}.$ In
particular $\frac{\partial u}{\partial y}(x,0)=0$ for $x\in
(-a,0).$ The difficulty in this approach is to prove the positivity of the
solution.
\end{remark}

\bigskip

\subsection{On the construction of some strict positive subsolutions with
zero flux on $(-a,0)\times \{0\}$}

The key point in the proof of Theorem \ref{Theorem No main} is the
construction of a flat positive subsolution. We start by reducing
the difficulty of this task by splitting the domain into two different
subdomains. We define the subsets 
\begin{equation*}
\Omega _{-}\text{ := }\left\{ (x,y)\in \Omega \text{ : }x\in (-a,0)\right\} ,%
\text{ and }\Omega _{+}\text{ := }\left\{ (x,y)\in \Omega \text{ : }x\in
(0,b)\right\} .
\end{equation*}

We will need to introduce an \textit{artificial boundary condition} on the
points $(0,y),$ $y\in (0,1)$, corresponding to the external boundary of this half of the cathode.
We will prove a rigorous result that, in some sense, is connected to the study made in \cite{Rokhlenko} by using asymptotic techniques and numerical analysis. Obtaining the correct
profile function in the external boundary condition, at $x=0$, is already an important conclusion,
since several choices were already considered in many papers in the literature on space charge problems. It was expected that the correct behaviour at the external border of the cathode is
greater than the profile $y^{4/3}$ in the middle of the cathode, corresponding to the Child-Langmuir law, but no concrete function has been proposed for this external profile. The precise artificial boundary condition on the exterior of the cathode, which we will introduce (as a Dirichlet boundary condition), is the following:

\begin{equation}
u(0,y)=hy^{\alpha },\text{ }y\in (0,1),\text{for some }h\in (0,1]\text{ and
for some }\alpha \in (1,4/3).  \label{Dirichlet Boundary condition}
\end{equation}

The following result will simplify the task of finding a global positive
subsolution since it will allow us to pass from the study of a discontinuous
absorption coefficient $j(x)$ to a problem with a strictly positive one.

\begin{proposition}
\label{Propo Matching}Let $j(x)$ satisfying (\ref{Hypo j a b}). Let $\alpha
\in (1,4/3).$ Given $h\in (0,1]$, consider the problem on $\Omega _{-}$ 
\begin{equation}
P_{a,0,j}=\left\{ 
\begin{array}{lr}
-\triangle u_{-}+\displaystyle\frac{j(x)}{\sqrt{u_{-}}}=0 & \text{in }\Omega
_{-}, \\ [.3cm]
u_{-}(-a,y)=y^{4/3} & y\in (0,1), \\ [.15cm]
u_{-}(0,y)=hy^{\alpha } & y\in (0,1), \\ [.15cm]
u_{-}(x,0)=0 & x\in (-a,0), \\ [.15cm]
u_{-}(x,1)=1 & x\in (-a,0).%
\end{array}%
\right.   \label{Box}
\end{equation}%
Assume that there exists $u_{0,-}(x,y)$, \textit{a }$p_{0}-$\textit{%
subsolution} of problem $P_{a,0,j}$, for some $p_{0}\geq 1$, such that 
\begin{equation}
\frac{\partial u_{0,-}}{\partial x}(0,y)\leq 0\text{ for }y\in (0,1),
\label{Hypo Subs decresing in Gamma}
\end{equation}%
satisfying the additional conditions 
\begin{equation}
AC_{a,0}=\left\{ 
\begin{array}{lr}
\frac{\partial u_{0,-}}{\partial y}(x,0)=0 & x\in (-a,0)\text{, } \\ [.15cm]
u_{0,-}(x,y)>0 & x\in (-a,0),\text{ }y\in (0,1).%
\end{array}%
\right.   \label{Additional condi Omega-}
\end{equation}%
Then problem $P_{a,b,j}$ has a $p_{0}-$\textit{subsolution} $\underline{u}$
satisfying the additional conditions $AC_{a,b}$.
\end{proposition}

\noindent Proof. Let $u_{+}(x,y)$ be the unique classical solution $u_{+}\in
C^{2}(\Omega _{+})\cap C^{0}(\overline{\Omega _{+}})$ of the linear problem 
\begin{equation}
P_{0,b,0}=\left\{ 
\begin{array}{lr}
-\triangle u_{+}=0\quad \text{ } & \text{in }\Omega _{+}, \\ [.3cm]
u_{+}(0,y)=hy^{\alpha }\quad & y\in (0,1), \\ [.15cm]
u_{+}(b,y)=y & y\in (0,1), \\ [.15cm]
u_{+}(x,0)=0 & x\in (0,b), \\ [.15cm]
u_{+}(x,1)=1 & x\in (0,b).%
\end{array}%
\right.  \label{Pb Omega+}
\end{equation}%
Define the function 
\begin{equation*}
\underline{u}(x,y)=\left\{ 
\begin{array}{lr}
u_{-}(x,y) & \text{if }(x,y)\in \Omega _{-}, \\ [.15cm]
u_{+}(x,y) & \text{if }(x,y)\in \Omega _{+}.%
\end{array}%
\right.
\end{equation*}

\noindent It is clear that $\underline{u}$ is a continuous function $%
\underline{u}\in C^{0}(\Omega )$ but its gradient has a discontinuity in the
segment $x=0,y\in (0,1)$, since $j(x)$ is discontinuous in that segment. All
the boundary conditions of $P_{a,b,j}$ are fulfilled and also the additional
conditions $AC_{a,b}$. In order to check that $\underline{u}(x,y)$ is a $%
p_{0}-$subsolution of problem $P_{a,0,j}$ we will apply Corollary I.1 of \ 
\cite{Berest-Lions}. To do this, if we define the segment $\Gamma =\{(0,y)$, 
$y\in (0,1)\}$ then it suffices to check that 
\begin{equation}
\displaystyle\frac{\partial u_{-}}{\partial \mathbf{n}}\leq \frac{\partial
u_{+}}{\partial \mathbf{n}}\quad \text{ on }\Gamma \text{,}
\label{dobleestrella}
\end{equation}%
where $\mathbf{n}$ is the unit exterior normal vector to $\Omega _{-}$. In
our case, $\mathbf{n}=\mathbf{e_{1}}$ and then condition (\ref{dobleestrella}%
) is expressed as 
\begin{equation}
\displaystyle\frac{\partial u_{-}}{\partial x}(0,y)\leq \frac{\partial u_{+}%
}{\partial x}(0,y)\qquad y\in (0,1).  \label{tripleestrella}
\end{equation}%
By the assumption (\ref{Hypo Subs decresing in Gamma}), it suffices to check that 
\begin{equation*}
\frac{\partial u_{+}}{\partial x}(0,y)\geq 0\text{ for }y\in (0,1).
\end{equation*}%
To do that, let us consider the function $U_{+}(x,y)=u_{+}(x,y)-y.$ Then%
\begin{equation}
\left\{ 
\begin{array}{lr}
-\triangle U_{+}=0\quad \text{ } & \text{in }\Omega _{+}, \\ [.3cm]
U_{+}(0,y)=y^{\alpha }-y & y\in (0,1), \\ [.15cm]
U_{+}(b,y)=0 & y\in (0,1), \\ [.15cm]
U_{+}(x,0)=0 & x\in (0,b), \\ [.15cm]
U_{+}(x,1)=0 & x\in (0,b).%
\end{array}%
\right.  \label{Problem U+}
\end{equation}

\noindent We can define, now, the auxiliary function,%
\begin{equation*}
\underline{U}(x,y)=(b-x)(hy^{\alpha }-y)\text{, }(x,y)\in \Omega _{+}.
\end{equation*}%
Then we get that 
\begin{equation*}
\left\{ 
\begin{array}{lr}
-\triangle \underline{U}\leq 0\quad \text{ } & \text{in }\Omega _{+}, \\ [.3cm]
\underline{U}(0,y)=\omega (y)=hy^{\alpha }-y & y\in (0,1), \\ [.15cm]
\underline{U}(b,y)=0 & y\in (0,1), \\ [.15cm]
\underline{U}(x,0)=0 & x\in (0,b), \\ [.15cm]
\underline{U}(x,1)=0 & x\in (0,b).%
\end{array}%
\right.
\end{equation*}

\noindent Thus, by the maximum principle we get that $\underline{U}(x,y)\leq
U_{+}(x,y)$ on $\Omega _{+}.$ But since $\underline{U}(0,y)=U_{+}(0,y)$ and
we have that $\frac{\partial \underline{U}}{\partial x}(0,y)\geq 0$ for $%
y\in (0,1),$ we deduce that necessarily, $\frac{\partial U_{+}}{\partial x}%
(0,y)\geq 0$ for $y\in (0,1),$ which leads to the required inequality.$\fin$

\bigskip

According Proposition \ref{Propo Matching}, to finish with the
construction of the global subsolution $u_{0}$, we must justify the
existence of a function $u_{0,-}(x,y)$ solution of the nonlinear problem $%
P_{a,0,j},$ raised on $\Omega _{-},$ and to check that $u_{0,-}(x,y)$ \
satisfies the additional conditions (\ref{Additional condi Omega-}) and (\ref%
{Hypo Subs decresing in Gamma}).

In order to construct the subsolution for the case of a possible unbounded $%
j(x),$ satisfying that $j(x)\leq \frac{A}{(-x)^{\beta }}$ on $(-a,0),$ we
will use some ideas coming from the study of Fluid Mechanics in the
consideration of spatial domains with corners (see, e.g., \cite{Batchelor}).
We will try to find the subsolution $u_{0,-}(x,y)$ in the form%
\begin{equation}
u_{0,-}(x,y)=\phi (r,\theta )=kr^{\alpha }U(\theta )\text{, for some }k>0%
\text{, }\alpha >1,  \label{subsolPower}
\end{equation}%
where $r=\sqrt{x^{2}+y^{2}}$ and $\theta =\arctan (y/x).$ Then, the partial
differential inequality becomes
\begin{equation}
-\triangle \phi +\frac{j(r\cos \theta )}{\sqrt{\phi }}=-\frac{1}{r}\frac{%
\partial }{\partial r}(r\frac{\partial \phi }{\partial r})-\frac{1}{r^{2}}%
\frac{\partial ^{2}\phi }{\partial \theta ^{2}}+\frac{j(r\cos \theta )}{%
\sqrt{\phi }}\leq 0.  \label{Ec EDP polares}
\end{equation}%
We will assume $r\in \lbrack 0,R)$, for a suitable $R>0$, and $\theta \in (%
\frac{\pi }{2},\pi ).$ The additional conditions (\ref{Additional condi
Omega-}) will require $\phi (r,\theta )>0$ if $r>0$ and 
\begin{equation*}
\phi (r,\pi )=\frac{\partial \phi }{\partial \theta }(r,\pi )=0\text{ for }%
r\in \lbrack 0,R).
\end{equation*}

\begin{figure}[htp]
\begin{center}
\includegraphics[width=12cm]{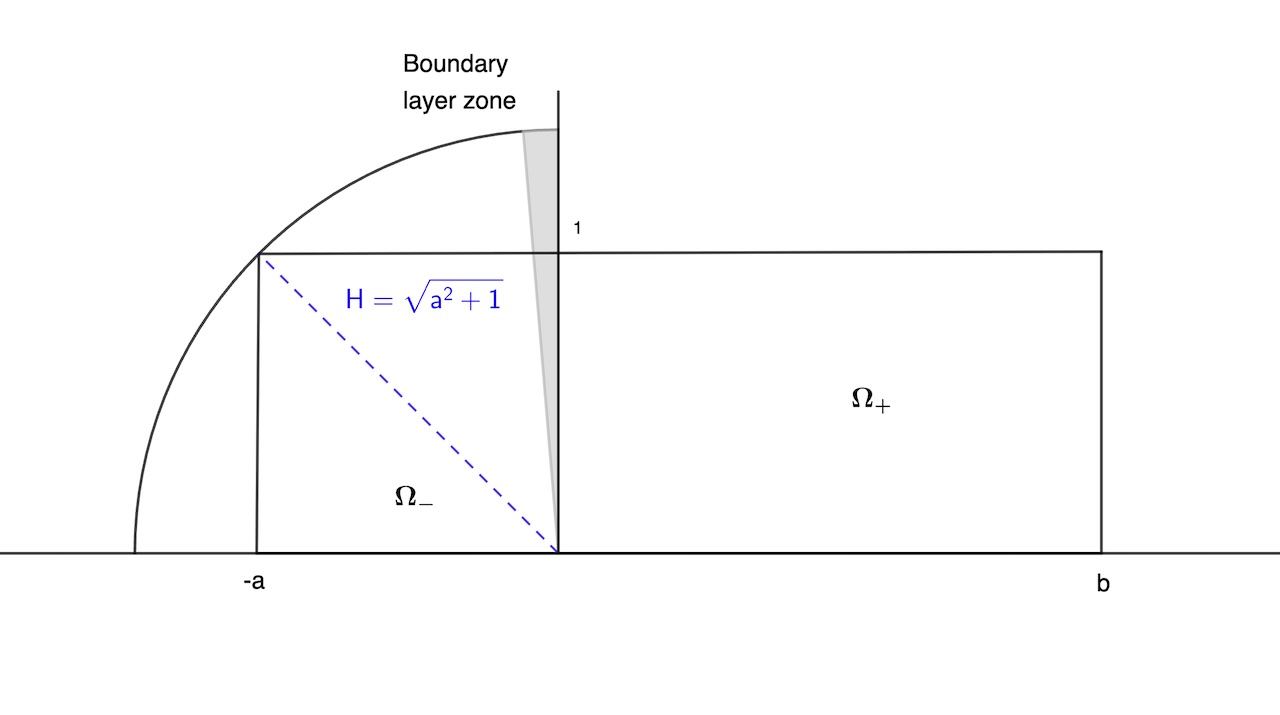}\\ 
\caption{Boundary layer}
\label{Fig boundary layer}
\end{center}
\end{figure}

We make now the \textit{structural condition} (\ref{Hypo j a b})\textbf{\ }%
on $j(x)$ and thus %
\begin{equation}
j(r\cos \theta )\leq \frac{A}{(-r\cos \theta )^{\beta }}\text{, if }\theta
\in (\frac{\pi }{2},\pi )\text{, for some }A>0\text{ and }\beta \in (0,1).
\label{Hypo j polar}
\end{equation}%
The partial differential inequality (\ref{Ec EDP polares}) leads to the
study of the ordinary differential equation

\begin{equation}
\begin{array}{lr}
-U^{\prime \prime }(\theta )+\frac{V(\theta )}{\sqrt{U(\theta )}}=\lambda
U(\theta ) & \theta \in (\frac{\pi }{2},\pi ),%
\end{array}
\label{Ec ODE3}
\end{equation}%
once we assume the constraint

\begin{equation}
\alpha +2\beta =4/3,  \label{Hypo alpha beta}
\end{equation}%
and then with 
\begin{equation}
V(\theta )=\frac{A}{(-\cos \theta )^{\beta }}\text{ and }\lambda =\alpha
^{2}.  \label{Hypo Potential and landa}
\end{equation}%
The complementary conditions $AC_{a,0}$ become now%
\begin{equation*}
\left\{ 
\begin{array}{cc}
U(\theta )>0 & \text{ }\theta \in (\frac{\pi }{2},\pi )\text{,} \\ [.15cm]
U(\pi )=U^{\prime }(\pi )=0. & 
\end{array}%
\right. 
\end{equation*}%
Notice that the potential $V(\theta )$ is singular only for $\theta =\frac{%
\pi }{2}$. On the other hand, we will need to match this subsolution with
another function which is positive for $\theta =\frac{\pi }{2}.$ Then, we will
construct the subsolution in two different pieces (see Figures \ref{Fig boundary layer} and \ref{fig Matching subsolution})%
\begin{equation*}
U(\theta )=\left\{ 
\begin{array}{cc}
v_{1}(\theta ) & \text{if }\theta \in (\frac{\pi }{2},\pi -R_{0})\text{,} \\ [.15cm]
v_{2}(\theta ) & \text{if }\theta \in (\pi -R_{0},\pi )\text{,}%
\end{array}%
\right. 
\end{equation*}%
for some $R_{0}\in (0,\frac{\pi }{2})$, and with $%
v_{2}(\theta )\in \lbrack 0,1]$ such that 
\begin{equation}
\left\{ 
\begin{array}{lc}
-v_{2}^{\prime \prime }(\theta )+\frac{V_{0}}{\sqrt{v_{2}(\theta )}}=\lambda
v_{2}(\theta ) & \theta \in (\pi -R_{0},\pi ), \\ 
v_{2}(\pi )=v_{2}^{\prime }(\pi )=0, & 
\end{array}%
\right.   \label{Equ ODE singular-constant}
\end{equation}%
where 
\begin{equation*}
V(\theta )\geq V_{0}\text{ if }\theta \in (\frac{\pi }{2},\pi -R_{0}).
\end{equation*}

On the other hand, $v_{1}(\theta )$ must take into account the singularity of
the potential $V(\theta )$ on the interval $(\frac{\pi }{2},\pi -R_{0})$ and
this will require a careful matching (see Figure \ref{Fig boundary layer}%
). In addition, we must guarantee a good match with the function $%
u_{+}(x,y)$ defined on $\Omega _{+}$, as indicated in Proposition %
\ref{Propo Matching}. This means that we want to have 
\begin{equation*}
v_{1}^{\prime }(\frac{\pi }{2})\geq 0,
\end{equation*}%
since we require $\frac{\partial u_{0,-}}{\partial x}(0,y)\leq 0$ for $y\in
(0,1),$ and from $u_{0,-}(0,y)=\left. kr^{\alpha }U(\theta )\right\vert
_{\theta =\frac{\pi }{2}}$, we have %
\begin{equation*}
\frac{\partial u_{0,-}}{\partial x}(0,y)=\cdots=\left. -Ar^{\alpha -1}U^{\prime
}(\theta )\right\vert _{\theta =\frac{\pi }{2}}\leq 0.
\end{equation*}

Before presenting the details on the construction of $v_{1}(\theta )$ and $%
v_{2}(\theta )$, it is very useful to consider the auxiliary nonlinear ODE of
eigenvalue type (\ref{Problem Psi}) presented in the Introduction.

\subsubsection{Bifurcation curve and flat solution for an auxiliary related
nonlinear eigenvalue ordinary differential problem}

Since the equation (\ref{Equ ODE singular-constant}) can be understood as a
nonlinear eigenvalue problem, it is useful to start by considering the
following auxiliary related problem:%
\begin{equation}
\left\{ 
\begin{array}{lr}
-U^{\prime \prime }(s)+\frac{V_{0}}{\sqrt{U(s)}}=\lambda U(s) & s\in (-R,R),
\\ 
U(\pm R)=0, & 
\end{array}%
\right.   \label{Problem nonlinear-eigenvalue}
\end{equation}%
where the positive constants $V_{0}$ and $R$ are given. The following result
extends several results in the previous literature (see, e.g., \cite%
{Diaz-Hernandez}, \cite{Di-Hern-Mance}, \cite{DHPortugalia} and some of the
expositions made in \cite{Brezis-Diaz-eds}).

\begin{theorem}
\label{Thm bifucation ODE} Given the positive constants $V_{0}$ and $R$ then:

i) there is a bifurcation from infinity for $\lambda $ near $\lambda
_{1}(R)=(\frac{\pi }{2R})^{2}$ (the first eigenvalue of the linear problem
with $V_{0}=0$),

ii) the bifurcation curve is strictly decreasing (which implies the
uniqueness of nonnegative solutions),

iii) the curve is not $C^{1}$ for a suitable value $\lambda =\lambda ^{\ast
}>\lambda _{1}(R)$ corresponding to a \textquotedblleft flat
solution\textquotedblright\ (i.e. the solution $U$ is such that $U^{\prime
}(\pm R)=0$ and $U(s)>0$).
\end{theorem}

To show the qualitative behaviour of solutions of
problem (\ref{Problem nonlinear-eigenvalue}) we make the change of variables 
\begin{equation*}
u_{\lambda ,V_{0}}(x)=\left( \frac{V_{0}}{\lambda }\right) ^{\frac{2}{3}}u( 
\sqrt{\lambda }x),
\end{equation*}%
where $u$ is now the solution of the renormalized problem 
\begin{equation}
P(L)\left\{ 
\begin{array}{lr}
-u^{\prime \prime }=f(u)\text{ } & \text{in }(-L,L), \\ 
u(\pm L)=0, & 
\end{array}
\right.  \label{problem L}
\end{equation}%
with
\begin{equation*}
f(u)=u-\frac{1}{\sqrt{u}}
\end{equation*}%
and $L=\sqrt{\lambda }R$. By multiplying by $u$ and integrating by parts, we
get that nontrivial solutions may exist only if $\lambda >\lambda _{1}=%
\frac{\pi ^{2}}{4R^{2}},$ the first eigenvalue to the linear problem 
\begin{equation}
\left\{ 
\begin{array}{lr}
-u^{\prime \prime }=\lambda u & \text{in }(-R,R), \\ 
u(\pm R)=0. & 
\end{array}
\right.  \label{linear problem}
\end{equation}

\begin{figure}[htp]
\begin{center}
\includegraphics[width=8cm]{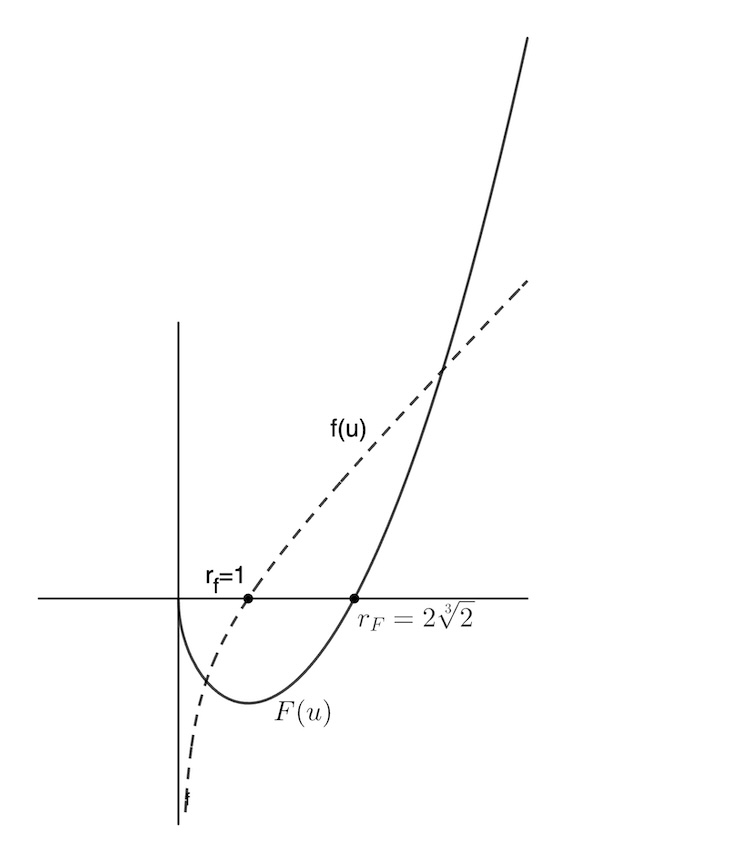}\\ 
\caption{Functions $f$ and $F.$}
\label{Figure f}
\end{center}
\end{figure}

\noindent Notice that now the role of the \textquotedblleft
eigenvalue\textquotedblright\ $\lambda $ is transferred to the length of the
new interval $L=\sqrt{\lambda }R.$ We introduce 
\begin{equation*}
F(r)=\int_{0}^{r}f(s)ds=\frac{r^{2}}{2}-2\sqrt{r},
\end{equation*}%
and note that $f(s)<0$ if $0<s<1:=r_{f}$ and $f(s)>0$ if $1<s$. On the other
hand $F(s)<0$ if $0<s<r_{F}=2\sqrt[3]{2}$ $(\approx 2,51)$ and $F(s)>0$ for $%
s>r_{F}$ (see Figure \ref{Figure f}).

\noindent By multiplying by $u^{\prime }$, integrating by parts and
denoting $\mu :=\left \Vert u\right \Vert _{L^{\infty }}$ for $\mu \in
[r_{F},\infty )$ we get that a function $u$\ is a positive solution of
problem $P(L)$\ if and only if%
\begin{equation*}
\frac{1}{\sqrt{2}}\int_{u(x)}^{\mu }\frac{dr}{(F(\mu )-F(r))^{1/2}}=\left
\vert x\right \vert ,\mbox{ for }\left \vert x\right \vert \leq {L},
\end{equation*}%
and $\mu $ and $L>0$\ are related by the equation 
\begin{equation*}
\gamma (\mu )=L,
\end{equation*}

\noindent where $\gamma :[r_{F},+\infty )\rightarrow \mathbb{R}$ is given by

\begin{equation}
\gamma (\mu ):=\frac{1}{\sqrt{2}}\int_{0}^{\mu }\frac{dr}{(F(\mu
)-F(r))^{1/2}}.
\end{equation}%


\noindent Now we use the following fact, whose proof is exactly as in \cite%
{Diaz-Hernandez} and \cite{Di-Hern-Mance}: a function $u$\ is a positive
solution of problem $P(L)$\ if and only if 
\begin{equation*}
\frac{1}{\sqrt{2}}\int_{u(x)}^{\mu }\frac{dr}{(F(\mu )-F(r))^{1/2}}%
=\left\vert x\right\vert ,\mbox{ for }\left\vert x\right\vert \leq {L}.
\end{equation*}%
Moreover

\begin{equation}
u^{\prime }(\pm L)=\mp \sqrt{2}\sqrt{F(\mu )}.  \label{derivada en L}
\end{equation}

\noindent Thus, $u^{\prime }(\pm R)=0$ corresponds to the case in which the
maximum of the solution is $r_{F}.$ It can be computed (using the \textit{\
Gauss-Lobatto rules: }a remark made by G. D\'{\i}az to the author) that $%
\gamma (r_{F})\approx 2,09$.

We will start by proving the existence of a branch of positive solutions for
a bounded interval of the parameter, $\lambda \in ((\frac{\pi }{2R}%
)^{2},\lambda _{1}^{\ast })$.

\begin{theorem}
\label{Thm bifurcation} \textit{We define with }$F(r)=\frac{r^{2}}{2}-2\sqrt{%
r}$\textit{. Let }$r_{F}=2\sqrt[3]{2}$\textit{. Then the mapping }$\gamma :[r_{F},+\infty )\rightarrow \mathbb{R}$ \textit{\ has the following properties}

\noindent \emph{(i) }$\gamma \in C[r_{F},\infty )\cap C^{1}(r_{F},\infty );$ 
\emph{\ }

\noindent \emph{(ii) For any }$\mu >r_{F}$ $\gamma ^{\prime }(\mu )<0$,

\noindent \emph{(iii) }$\gamma ^{\prime }(\mu )\rightarrow -\xi $\emph{\ as }
$\mu \downarrow r_{F}$\emph{, for some }$\xi >0$\emph{, }

\noindent \emph{(iv) }$\lim_{\mu \rightarrow +\infty }\gamma (\mu )=\frac{\pi }{2}.$
\end{theorem}

\bigskip Concerning the flat solution and the possible solutions with compact support we have:

\begin{theorem}
\label{Thm bifurcation2}
\noindent \textit{Let } 
\begin{equation}
\lambda _{1}^{\ast }=\frac{1}{2R^{2}}\left( \int_{0}^{r_{F}}\frac{dr}{(F(\mu
)-F(r))^{1/2}}\right) ^{2}  \label{landa critico}
\end{equation}%
\textit{then we have:}

\textit{\noindent a) if }$\lambda \in (0,(\frac{\pi }{2R})^{2})$\textit{\
there is no positive solution,}

\textit{\noindent b) if }$\lambda \in ((\frac{\pi }{2R})^{2},\lambda
_{1}^{\ast })$\textit{\ there is a unique positive solution }$u_{\lambda
,V_{0}}$\textit{. Moreover }$\partial u_{\lambda ,V_{0}}/\partial n(\pm R)<0$
\textit{\ and }$\left\Vert u_{\lambda ,V_{0}}\right\Vert _{L^{\infty
}(-R.R)}=\left( \frac{V_{0}}{\lambda }\right) ^{\frac{2}{3}}\gamma ^{-1}(%
\sqrt{\lambda }R),$

\textit{\noindent c) if }$\lambda =\lambda _{1}^{\ast }$\textit{\ there is
only one positive solution }$\ u_{\lambda _{1}^{\ast },V_{0}}$\textit{.
Moreover }$u_{\lambda _{1}^{\ast },V_{0}}^{\prime }(\pm R)=0$\textit{\ } 
\begin{equation*}
\left\Vert u_{\lambda _{1}^{\ast },V_{0}}\right\Vert _{L^{\infty
}(-R,R)}=\left( \frac{4V_{0}}{\lambda _{1}^{\ast }}\right) ^{\frac{2}{3}}.
\end{equation*}

\textit{\noindent d) if }$\lambda >\lambda _{1}^{\ast },$\textit{\ there is
a family of nonnegative solutions that are generated by extending by zero
the function }$u_{\lambda _{1}^{\ast },V_{0}}$\textit{\ outside }$(-R,R)$ 
\textit{\ (and which we label again as }$u_{\lambda _{1}^{\ast },V_{0}}$). In particular, if $\lambda =\lambda _{1}^{\ast }\omega $\textit{%
\ \ with }$\omega >1$\textit{\ we have a family }$S_{1}(\lambda )$\textit{\
of compact support nonnegative solutions with connected support defined by } 
\begin{equation*}
u_{\lambda ,V_{0}}(x)=\frac{1}{\omega ^{^{\frac{2}{3}}}}u_{\lambda
_{1}^{\ast },V_{0}}(\sqrt{\omega }x-z),
\end{equation*}%
\textit{where the shifting argument }$z$\textit{\ is arbitrary among the
points }$z\in (-R,R)$\textit{\ such that support of }$u_{\lambda
,V_{0}}(.)\subset (-R,R).$\textit{\ Moreover, for }$\lambda >\lambda
_{1}^{\ast }$\textit{\ large enough we can build, similarly, a subset of }$%
S_{j}(\lambda )$\textit{\ of compact support nonnegative solutions with the
support formed by }$j$components\textit{, with }$j\in \{1,2,...,N\},$\textit{%
\ \ for some suitable }$N=N(\lambda )$\textit{\ and then the set of
nontrivial and nonnegative solutions of }$P(\lambda )$\textit{\ is formed by 
}$S(\lambda )=\cup _{j=1}^{N}S_{j}(\lambda ).$\textit{\ In any case those
solutions satisfy that } 
\begin{equation*}
\left\Vert u_{\lambda ,V_{0}}\right\Vert _{L^{\infty }(-R,R)}=\frac{1}{%
\omega ^{^{\frac{2}{3}}}}\left\Vert u_{\lambda _{1}^{\ast
},V_{0}}\right\Vert _{L^{\infty }(-R,R)}=\frac{1}{\omega ^{^{\frac{2}{3}}}}%
\left( \frac{4V_{0}}{\lambda _{1}^{\ast }}\right) ^{\frac{2}{3}}\text{, 
\textit{for any} }\omega =\lambda /\lambda _{1}^{\ast }>1.
\end{equation*}
\end{theorem}

\noindent \textit{Proofs of Theorems \ref{Thm bifurcation} and \ref{Thm bifurcation2}}. The proof of
property i) is exactly the same as the one presented in \cite%
{Di-Hern-Mance}. For the proof of (ii) and (iii), we have 
\begin{equation*}
\gamma ^{\prime }(\mu )=\int_{0}^{\mu }\frac{\theta (\mu )-\theta (r)}{
(F(\mu )-F(r))^{1/2}}dr
\end{equation*}%
where $\theta (t)=2F(t)-tf(t)=-3\sqrt{t},$ and differentiating we get for,
any $t>0$%
\begin{equation*}
\theta ^{\prime }(t)=-\frac{3}{2\sqrt{t}}<0.
\end{equation*}%
Hence $\gamma ^{\prime }(\mu )<0$ for any $\mu >r_{F}$ (which proves (ii)).

\noindent For the proof of (iii) it suffices to see that as $\mu \downarrow {%
\ r}_{F},$ the integrand of $\gamma ^{\prime }(\mu )$ converges pointwise to 
$(-F(r))^{-3/2}$ near $r=0$ and in our case $(-F(r))^{-3/2}$ behaves as $%
r^{-3/4}$ near $r=0$ and thus $\gamma ^{\prime }(\mu )$ converges to a
number $-\xi $\emph{\ }as $\mu \downarrow r_{F}$, for some $\xi >0$.

\noindent Finally, to prove (iv), we note that 
\begin{equation}
\gamma (\mu )\leq \frac{\mu }{2}\int_{0}^{1}\frac{dt}{(\frac{\mu ^{2}}{2}
(1-t^{2}))^{1/2}}=\int_{0}^{1}\frac{dt}{\sqrt{1-t^{2}}}=\frac{\pi }{2}.
\end{equation}%
Moreover, we have 
\begin{equation}
\gamma (\mu )=\frac{\mu }{\sqrt{2}}\int_{0}^{1}\frac{dt}{\left( \frac{\mu }{ 
\sqrt{2}}((1-t^{2})-\frac{1}{\mu ^{3/2}}(1-\sqrt{t^{{}}})\right) },
\end{equation}%
and if $\mu \rightarrow +\infty, $ by using Lebesgue's Theorem, we get 
\begin{equation*}
\lim_{\mu \rightarrow +\infty }\gamma (\mu )=\frac{\pi }{2}.
\end{equation*}

\noindent\ Now we define $L_{0}=\frac{\pi }{2}$ and $L^{\ast }$ given by 
\begin{equation*}
L^{\ast }=\gamma (r_{F})=\frac{1}{\sqrt{2}}\int_{0}^{r_{F}}\frac{dr}{
(-F(r))^{1/2}}=\frac{1}{\sqrt{2}}\int_{0}^{r_{F}}\frac{dr}{(2\sqrt{r}-\frac{
r^{2}}{2})^{1/2}}.
\end{equation*}%
We know that $u^{\prime }(\pm R)=0$ corresponds to the value $L^{\ast }$ and
that the maximum of the solution is $r_{F}.$

\begin{figure}[htp]
\begin{center}
\includegraphics[width=8cm]{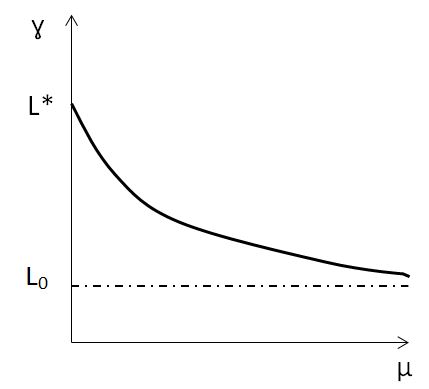}\\ 
\caption{Critical value $L^{\ast }.$}
\label{Fig gamma}
\end{center}
\end{figure}

\noindent So, qualitatively, function $\gamma $ is described by the Figure \ref{Fig gamma}.

\noindent If we now go back with our change of variables, we get 
\begin{equation*}
\left\Vert u_{\lambda ,V_{0}}\right\Vert _{L^{\infty }(-R,R)}=\left( \frac{%
V_{0}}{\lambda }\right) ^{\frac{2}{3}}r_{F},
\end{equation*}%
and we obtain, finally, the bifurcation diagram given by the first branch of
Figure \ref{Fig Canarias}, where solutions for $\lambda >\lambda ^{\ast }$
are compact supported solutions originated as in \cite{Di-Hern-Mance} from
the extension by zero of the free boundary solution $u_{\lambda ^{\ast }}$
satisfying 
\begin{equation}
\left\{ 
\begin{array}{lr}
-u_{\lambda ^{\ast }}^{\prime \prime }+\frac{V_{0}}{\sqrt{u_{\lambda ^{\ast
}}}}=\lambda ^{\ast }u_{\lambda ^{\ast }}\text{ } & \text{in }(-R,R), \\ [.25cm]
u_{\lambda ^{\ast }}(\pm R)=u_{\lambda ^{\ast }}^{\prime }(\pm R)=0. & 
\end{array}%
\right. 
\end{equation}%
The rest of the details are completely analogous to the similar parts of the paper 
\cite{Di-Hern-Mance}.$\fin$

\begin{remark}
\noindent Once that we know that for $\lambda >\lambda _{1}^{\ast }$ we have
that\textit{\ } 
\begin{equation*}
u_{\lambda ,V_{0}}(x)=\frac{1}{\omega ^{^{\frac{2}{3}}}}u_{\lambda
_{1}^{\ast },V_{0}}(\sqrt{\omega }x-z),
\end{equation*}
then we get that the bifurcating curve $\Lambda $ is not $C^{1}$ at $\lambda
=\lambda _{1}^{\ast }$ since $\Lambda ^{\prime }(\lambda _{1}^{\ast }-)=\xi
<0$ and $\Lambda ^{\prime }(\lambda _{1}^{\ast }+)=-\frac{2C}{3(\lambda
_{1}^{\ast })^{5/3}},$ with $C=4^{2/3}(V_{0})^{2/3}.$ In addition, for $%
\lambda >\lambda _{1}^{\ast }$ we can express other norms (different from
the $L^{\infty }$-norm) in terms of $\lambda $. For instance, we have that 
\begin{equation*}
\left\Vert u_{\lambda ,V_{0}}^{\prime }\right\Vert _{L^{\infty
}(-R,R)}=C\lambda ^{-\frac{1}{6}}
\end{equation*}
for a suitable constant $C>0$ independent of $\lambda $. This proves that $%
\left\Vert u_{\lambda ,V_{0}}^{\prime }\right\Vert _{L^{\infty
}(-R,R)}\rightarrow 0$ as $\lambda \rightarrow +\infty $.
\end{remark}

\bigskip

\subsubsection{On the construction of the subsolution over $\Omega _{-}$:
continuation}
\begin{figure}[htp]
\begin{center}
\includegraphics[width=8cm]{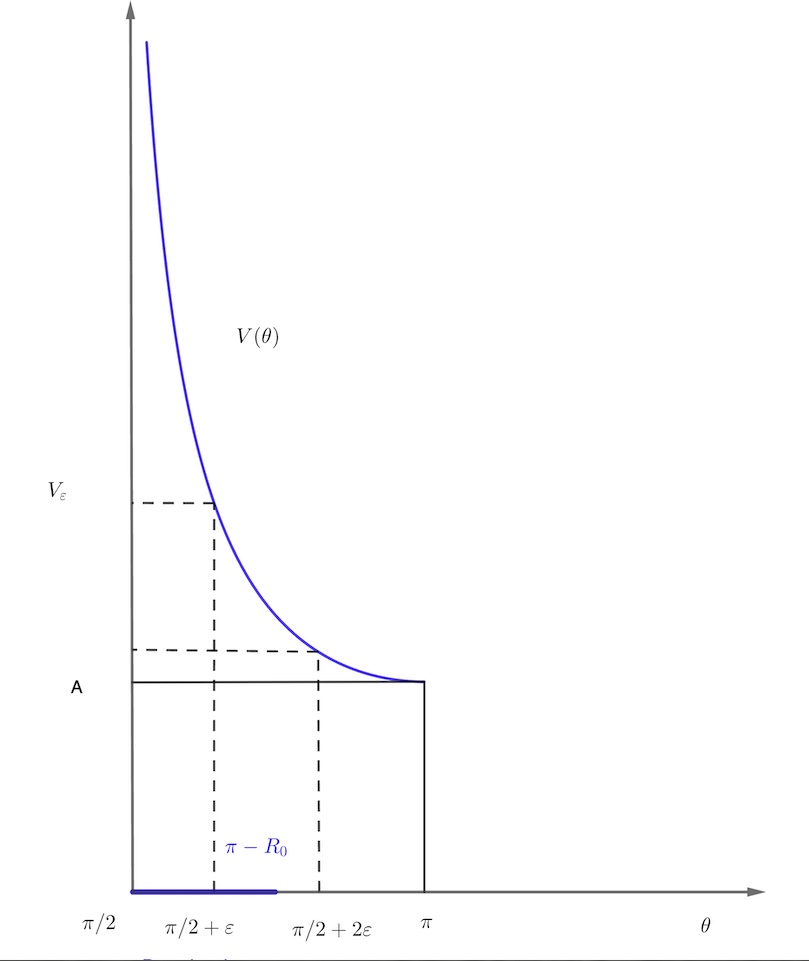}\\ 
\caption{Singular potential function $V(\theta ).$ }
\label{Figure V singular}
\end{center}
\end{figure}

\textit{Proof of Theorem 27}. Now, let us come back to the construction of the positive subsolution over $%
\Omega _{-}$ mentioned in Proposition \ref{Propo Matching}.\ We recall that,
for functions of the form (\ref{subsolPower}) the partial differential
equation (\ref{Ec EDP polares}) leads to the ordinary differential equation

\begin{equation}
\begin{array}{lr}
-U^{\prime \prime }(\theta )+\frac{V(\theta )}{\sqrt{U(\theta )}}=\lambda
U(\theta ) & \theta \in (\frac{\pi }{2},\pi ),%
\end{array}
\label{Ec ODE}
\end{equation}%
once we assume the constraint

\begin{equation}
\alpha +2\beta =4/3,
\end{equation}%
with 
\begin{equation}
V(\theta )=\frac{A}{(-\cos \theta )^{\beta }}\text{ and }\lambda =\alpha
^{2}.
\end{equation}%
Notice that the potential $V(\theta )$ is singular only for $\theta =\frac{%
\pi }{2}$. On the other hand, we need to match this subsolution with another
function which is positive for $\theta =\frac{\pi }{2}$. Then, given $%
\varepsilon >0$ small enough, we will construct the subsolution in two
different pieces%
\[
U(\theta )=\left \{ 
\begin{array}{cc}
v_{1}(\theta ) & \text{if }\theta \in (\frac{\pi }{2},\pi -R_{0})\text{,} \\%
[.15cm] 
v_{2}(\theta ) & \text{if }\theta \in (\pi -R_{0},\pi )\text{,}%
\end{array}%
\right.  \label{Matched subsoluion}
\]%
for some $R_{0}$ with $\pi -R_{0}\in (\frac{\pi }{2}+\varepsilon ,\frac{\pi 
}{2}+2\varepsilon )$, and with $v_{2}(\theta )$ such that 
\begin{equation}
\left \{ 
\begin{array}{lc}
-v_{2}^{\prime \prime }(\theta )+\frac{V_{\varepsilon }}{\sqrt{v_{2}(\theta )%
}}=\lambda v_{2}(\theta ) & \theta \in (\pi -R_{0},\pi ), \\[.3cm] 
v_{2}(\pi )=v_{2}^{\prime }(\pi )=0, & 
\end{array}%
\right.
\end{equation}%
where 
\[
V(\theta )\leq V_{\varepsilon }\text{ if }\theta \in (\frac{\pi }{2}%
+\varepsilon ,\pi ).
\]%
\noindent See Figure \ref{Figure V singular}.

\noindent \textit{Proof of Theorem \ref{Theorem No main}}. We recall the
change of variable 
\[
u_{\lambda ,V_{0}}(x)=\left( \frac{V_{0}}{\lambda }\right) ^{\frac{2}{3}}u(%
\sqrt{\lambda }x),
\]%
which links the question of the distinguished eigenvalue $\lambda ^{\ast
}=(\alpha ^{\ast })^{2}=\frac{1}{2R^{2}}\gamma (r_{F})^{2}$ with a special
length $L^{\ast }$ for which the solution is flat on the boundary.

\noindent The computations $r_{F}=4^{2/3}$ $(\approx 2,51)$ and $\gamma
(r_{F})\approx 2,09 $ allow to see that the corresponding $L^{\ast }$ leads
to $\lambda ^{\ast }=\frac{16}{9}$ and thus $\alpha ^{\ast }=\frac{4}{3}$
and $R=\frac{\pi }{2}$. \ This corresponds to the case $\beta =0$ thanks to
the reciprocal relation 
\[
\beta =2-\frac{3\alpha }{2}.
\]%
Moreover, we must take $A$ small enough (in fact $A=\frac{4}{9}$) and there
is a kind of \textit{fortunate coincidence}.

\bigskip

\noindent Notice that although $\left \Vert u_{\lambda ,V_{0}}\right \Vert
_{L^{\infty }}>1$ there is no difficulty with this since we can take $%
u_{0,-}(x,y)=\phi (r,\theta )=kr^{\alpha }U(\theta )$, for some $k>0$ small
enough.

\begin{figure}[htp]
\begin{center}
\includegraphics[width=10cm]{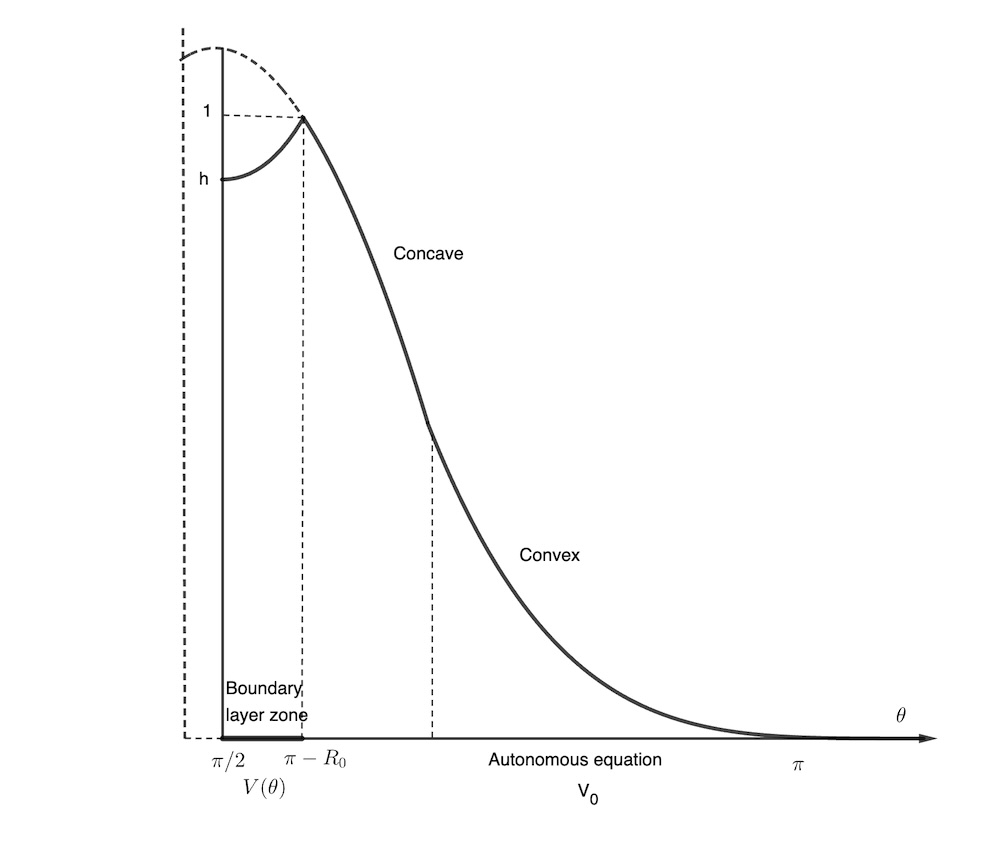}\\[0pt]
\end{center}
\caption{Matching functions $v_{1}(\protect \theta )$ and $v_{2}(\protect%
\theta )$ to get a subsolution.}
\label{fig Matching subsolution}
\end{figure}

\noindent If we consider the case of $\beta \in (0,\frac{1}{2})$ (i.e. $%
\alpha ^{\ast }\in (1,\frac{4}{3})$) then the problem is not autonomous (it
appears $V(\theta )=\frac{A}{(-\cos \theta )^{\beta }}$), the corresponding $%
R>0$, given by the equation $\lambda ^{\ast }=(\alpha ^{\ast })^{2}=\frac{1}{%
2R^{2}}\gamma (r_{F})^{2}$, is such that $R>\frac{\pi }{2}$ and we must
truncate $V(\theta )$ on an interval $(\pi -R_{0},\pi )$, for instance by
taking $R_{0}<R$ such that $u_{\lambda ,V_{0}}(\pi -R_{0})=1$ (see Figures %
\ref{Figure V singular} and \ref{fig Matching subsolution}).

\noindent We see then that the matching between $v_{1}(\theta )$ and $%
v_{2}(\theta )$ must take place at $\pi -R_{0}$ (see Figure \ref{fig
Matching subsolution}).

\bigskip

\noindent The construction of $v_{1}(\theta )$, in the so-called
\textquotedblleft boundary layer zone\textquotedblright , can be carried out
as in the one-dimensional problem with a distributed potential $V(\theta )$.
Indeed, given $h\in (0,1)$ we take $v_{1}(\theta )$ as the solution of%
\begin{equation}
\left \{ 
\begin{array}{lc}
-v_{1}^{\prime \prime }(\theta )+\frac{V(\theta )}{\sqrt{v_{1}(\theta )}}%
\leq 0 & \theta \in (\frac{\pi }{2},\pi -R_{0}), \\ 
v_{1}(\frac{\pi }{2})=h\text{, }v_{1}(\pi -R_{0})=K\in (h,1]. & 
\end{array}%
\right.
\end{equation}

\noindent Notice that using that%
\begin{equation*}
\cos \theta \leq \frac{2}{\pi }(\frac{\pi }{2}-\theta )\text{, if }\theta
\in (\frac{\pi }{2},\pi ),
\end{equation*}%
then the structural assumption on $V(\theta )$ implies that 
\begin{equation*}
V(\theta )\leq \frac{A\pi ^{\beta }}{2^{\beta }(\theta -\frac{\pi }{2}%
)^{\beta }}:=\frac{C}{(\theta -\frac{\pi }{2})^{\beta }}.
\end{equation*}%
and we can take%
\begin{equation*}
-v_{1}^{\prime \prime }(\theta )+\frac{\frac{C}{(\theta -\frac{\pi }{2}%
)^{\beta }}}{\sqrt{v_{1}(\theta )}}\leq 0
\end{equation*}%
with the above boundary conditions. Notice that this is an important key
point in the proof, since in Corollary 19, the possible singular coefficient
was depending on the $y$ variable but the coefficient $j$ in the
two-dimensional problem depends of the other variable $j=j(x).$ By applying
Corollary 19 to $q=-\beta \in (-\frac{1}{2},0)$, we can define%
\begin{equation*}
v_{1}(\theta )=u_{q}(\theta -\frac{\pi }{2})+h
\end{equation*}%
and thus, $v_{1}(\theta )\geq u_{q}(y),$ if $y=\theta -\frac{\pi }{2},$ 
\begin{equation*}
-v_{1}^{\prime \prime }(\theta )+\frac{V(\theta )}{\sqrt{v_{1}(\theta )}}%
\leq -v_{1}^{\prime \prime }(\theta )+\frac{\frac{C}{(\theta -\frac{\pi }{2}%
)^{\beta }}}{\sqrt{v_{1}(\theta )}}\leq -u_{q}^{\prime \prime }(y)+\frac{%
Cy^{q}}{\sqrt{u_{q}}}=0
\end{equation*}%
Thus, we have 
\begin{equation*}
v_{1}^{\prime }(\frac{\pi }{2})=0.
\end{equation*}%
We get that the graph of the matched function presents a cusp at $\theta
=\pi -R_{0}$, which is far from being \textit{natural}, but it is sufficient
for our purposes. Indeed, we can regularize the peaked function by means of
a third transition function $w(\theta )$ if $\theta \in (\frac{\pi }{2},%
\frac{\pi }{2}+2\varepsilon )$, so that the resulting subsolution will be of
the form 
\begin{equation}
U(\theta )=\left \{ 
\begin{array}{cc}
v_{1}(\theta ) & \text{if }\theta \in (\frac{\pi }{2},\frac{\pi }{2}%
+\varepsilon ]\text{,} \\ 
w(\theta ) & \text{if }\theta \in (\frac{\pi }{2}+\varepsilon ,\frac{\pi }{2}%
+2\varepsilon )\text{,} \\ 
v_{2}(\theta ) & \text{if }\theta \in \lbrack \frac{\pi }{2}+2\varepsilon
,\pi )\text{.}%
\end{array}%
\right.   \label{Subsolutin 3funct}
\end{equation}%
We start by recalling that the change of variables made in the proof of
Theorem 32 shows that the parameter $\lambda ^{\ast }=(\alpha ^{\ast })^{2}=%
\frac{1}{2R^{2}}\gamma (r_{F})^{2}$ arising in the problem definition of $%
v_{2}(\theta )$ is independent of the coefficient $V_{\varepsilon }$ also
arising in the nonlinear term of this problem. Thus, we will take $%
\varepsilon >0$ such that 
\begin{equation*}
V(\theta )\leq V_{\varepsilon }<\lambda ^{\ast },\text{ if }\theta \in (%
\frac{\pi }{2}+\varepsilon ,\pi ).
\end{equation*}%
We also recall that, from the above construction, we assume $\beta >0$ small
enough (which implies that also $R-\frac{\pi }{2}$ is small enough) and then
there exists $\delta >0$ such that, if $\theta \in (\frac{\pi }{2}%
+\varepsilon ,\frac{\pi }{2}+2\varepsilon ),$ 
\begin{eqnarray*}
\max (1-v_{1}(\theta ),1-v_{2}(\theta )) &\leq &\delta ,\text{ and } \\
\max (v_{1}^{\prime }(\theta ),-v_{2}^{\prime }(\theta )) &\leq &\delta .
\end{eqnarray*}%
We take $R_{0}=\frac{\pi }{2}+\frac{3\varepsilon }{2}.$ Then, we define the
concave function 
\begin{equation*}
w(\theta )=k_{0}+k_{1}(\theta -R_{0})-k_{2}(\theta -R_{0})^{2},\text{ if }%
\theta \in (\frac{\pi }{2}+\varepsilon ,\frac{\pi }{2}+2\varepsilon )=(R_{0}-%
\frac{\varepsilon }{2},R_{0}+\frac{\varepsilon }{2}),
\end{equation*}%
for some constants $k_{0},k_{1},k_{2}$ to be chosen now, with $k_{0},k_{2}>0$%
. On the interval $(\frac{\pi }{2}+\varepsilon ,\frac{\pi }{2}+2\varepsilon )
$ we have 
\begin{equation*}
w(\theta )\geq k_{0}-\left \vert k_{1}\right \vert \frac{\varepsilon }{2}-%
\frac{k_{2}\varepsilon ^{2}}{4}.
\end{equation*}%
Then, the subsolution condition, concerning the differential equation, leads
to the inequality%
\begin{equation}
-w^{\prime \prime }(\theta )+\frac{V_{\varepsilon }}{\sqrt{w(\theta )}}%
-\lambda ^{\ast }w(\theta )\leq 2k_{2}+\frac{V_{\varepsilon }}{\sqrt{%
k_{0}-\left \vert k_{1}\right \vert \frac{\varepsilon }{2}-\frac{%
k_{2}\varepsilon ^{2}}{4}}}-\lambda ^{\ast }(k_{0}-\left \vert
k_{1}\right \vert \frac{\varepsilon }{2}-\frac{k_{2}\varepsilon ^{2}}{4})\leq
0.  \label{Hypo wSubsolution}
\end{equation}%
On the boundary, we must demand the continuity of the matching%
\begin{equation}
\left \{ 
\begin{array}{cc}
w(\frac{\pi }{2}+\varepsilon )=k_{0}-k_{1}\frac{\varepsilon }{2}-\frac{%
k_{2}\varepsilon ^{2}}{4}=v_{1}(\frac{\pi }{2}+\varepsilon ), &  \\ 
w(\frac{\pi }{2}+2\varepsilon )=k_{0}+k_{1}\frac{\varepsilon }{2}-\frac{%
k_{2}\varepsilon ^{2}}{4}=v_{2}(\frac{\pi }{2}+2\varepsilon ), & 
\end{array}%
\right.   \label{eq. Continuity matching subsol}
\end{equation}%
and the good sign conditions on the measures generated by the derivatives
matching%
\begin{equation}
\left \{ 
\begin{array}{cc}
w^{\prime }(\frac{\pi }{2}+\varepsilon )=k_{1}+2k_{2}\frac{\varepsilon }{2}%
\geq v_{1}^{\prime }(\frac{\pi }{2}+\varepsilon ), &  \\ 
w^{\prime }(\frac{\pi }{2}+2\varepsilon )=k_{1}-2k_{2}\frac{\varepsilon }{2}%
\leq v_{2}^{\prime }(\frac{\pi }{2}+2\varepsilon ). & 
\end{array}%
\right.   \label{Eq- Matching deriv subsol}
\end{equation}%
Let 
\begin{equation*}
v_{1}^{L}:=v_{1}\! \left( \tfrac{\pi }{2}+\varepsilon \right) ,\qquad
v_{2}^{R}:=v_{2}\! \left( \tfrac{\pi }{2}+2\varepsilon \right) .
\end{equation*}%
From (\ref{eq. Continuity matching subsol}) we have 
\begin{equation*}
k_{0}-k_{1}\frac{\varepsilon }{2}-\frac{k_{2}\varepsilon ^{2}}{4}=v_{1}^{L},%
\text{ }k_{0}+k_{1}\frac{\varepsilon }{2}-\frac{k_{2}\varepsilon ^{2}}{4}%
=v_{2}^{R}.
\end{equation*}%
Adding and subtracting these equations gives 
\begin{equation*}
k_1 = \frac{v_2^{R}-v_1^{L}}{\varepsilon},\qquad k_0 =
\frac{v_1^{L}+v_2^{R}}{2} + \frac{k_2\varepsilon^2}{4}.
\end{equation*}%
Hence $k_{1}$ is fixed by the endpoint difference and $k_{0}$ depends
linearly on $k_{2}>0$. From (\ref{Eq- Matching deriv subsol}) we have 
\begin{equation*}
k_{1}+k_{2}\varepsilon \geq v_{1}^{\prime }(L),\newline
k_{1}-k_{2}\varepsilon \leq v_{2}^{\prime }(R).
\end{equation*}%
Rearranging gives the lower bounds 
\begin{equation*}
k_{2}\geq \frac{v_{1}^{\prime }(L)-k_{1}}{\varepsilon },\qquad k_{2}\geq 
\frac{k_{1}-v_{2}^{\prime }(R)}{\varepsilon }.
\end{equation*}%
Thus any admissible $k_{2}$ must satisfy 
\begin{equation*}
k_{2}\geq k_{2,\min }:=\max \! \left \{ 0,\, \frac{v_{1}^{\prime }(L)-k_{1}}{%
\varepsilon },\, \frac{k_{1}-v_{2}^{\prime }(R)}{\varepsilon }\right \} .
\end{equation*}%
Using the hypotheses $|v_{1}^{\prime }(L)|,|v_{2}^{\prime }(R)|\leq \delta $
and $|v_{2}^{R}-v_{1}^{L}|\leq \delta $, we have $|k_{1}|\leq \delta
/\varepsilon $, so $k_{2,\min }=O(\delta /\varepsilon )$, typically small.
By substituting the relations for $k_{0},k_{1}$ into \eqref{ineq:main} we
can define: 
\begin{equation*}
\Phi (k_{2}):=2k_{2}+\frac{V_{\varepsilon }}{\sqrt{k_{0}-|k_{1}|(\varepsilon
/2)-\frac{k_{2}\varepsilon ^{2}}{4}}}-\lambda ^{\ast }\! \left(
k_{0}-|k_{1}|(\varepsilon /2)-\frac{k_{2}\varepsilon ^{2}}{4}\right) .
\end{equation*}%
Note that 
\begin{equation*}
k_{0}-|k_{1}|(\varepsilon /2)-\frac{k_{2}\varepsilon ^{2}}{4}=\frac{%
v_{1}^{L}+v_{2}^{R}}{2}-|k_{1}|\frac{\varepsilon }{2},
\end{equation*}%
which is independent of $k_{2}$. Hence, the desired inequality becomes 
\begin{equation*}
\Phi (k_{2})=2k_{2}+C\leq 0,
\end{equation*}%
where 
\begin{equation*}
C:=\frac{V_{\varepsilon }}{\sqrt{\frac{v_{1}^{L}+v_{2}^{R}}{2}%
-|k_{1}|(\varepsilon /2)}}-\lambda ^{\ast }\! \left( \frac{v_{1}^{L}+v_{2}^{R}%
}{2}-|k_{1}|(\varepsilon /2)\right) .
\end{equation*}%
Using $\max (1-v_{1},1-v_{2})\leq \delta $, we have $v_{1}^{L},v_{2}^{R}\geq
1-\delta $ and thus 
\begin{equation*}
\frac{v_{1}^{L}+v_{2}^{R}}{2}-|k_{1}|\frac{\varepsilon }{2}\geq 1-\frac{%
3\delta }{2}>0
\end{equation*}%
for small $\delta $. Therefore 
\begin{equation*}
-C\approx V_{\varepsilon }+\lambda ^{\ast }<0.
\end{equation*}%
By continuity, $C$ remains positive for sufficiently small $\delta $.
Because $C<0$, the function $\Phi (k_{2})=2k_{2}+C$ is increasing and linear
in $k_{2}$. Thus the inequality $\Phi (k_{2})\leq 0$ holds for any 
\begin{equation*}
0<k_{2}\leq \frac{-C}{2}.
\end{equation*}%
To satisfy both the derivative and main inequalities, it suffices to have 
\begin{equation*}
k_{2,\min }\leq k_{2}\leq \frac{-C}{2}.
\end{equation*}%
Since $k_{2,\min }=O(\delta /\varepsilon )$ and $-C/2\approx (\lambda ^{\ast
}-V_{\varepsilon })/2>0$, this interval is nonempty for small enough $\delta 
$. In conclusion, defining 
\begin{equation*}
B:=\frac{v_{1}^{L}+v_{2}^{R}}{2}-|k_{1}|\frac{\varepsilon }{2},\qquad C:=%
\frac{V_{\varepsilon }}{\sqrt{B}}-\lambda ^{\ast }B,
\end{equation*}%
a convenient, explicit choice is 
\begin{equation*}
\left \{ 
\begin{array}{ll}
k_{1}=\dfrac{v_{2}^{R}-v_{1}^{L}}{\varepsilon }, & k_{2,\min }=\max
\! \left \{ 0,\, \dfrac{v_{1}^{\prime }(L)-k_{1}}{\varepsilon },\, \dfrac{%
k_{1}-v_{2}^{\prime }(R)}{\varepsilon }\right \} , \\ 
k_{2}=k_{2,\min }+\min \! \left \{ \eta ,\, \frac{-C}{4}\right \} , & k_{0}=%
\dfrac{v_{1}^{L}+v_{2}^{R}}{2}+\dfrac{k_{2}\varepsilon ^{2}}{4},%
\end{array}%
\right. 
\end{equation*}
where $\eta >0$ is any small tolerance (e.g.,\ $\eta =\tfrac{1}{10}%
\varepsilon ^{3}$). Then, $k_{0},k_{2}>0$, all the required inequalities are
satisfied and the function $U(\theta )$ defined in (\ref{Subsolutin 3funct})
is a subsolution as required, $0<U(\theta )\leq 1$ for $\theta \in (\frac{%
\pi }{2},\pi ).$
Finally, the boundary inequality
\begin{equation*}
u_{0,-}(x,y)=kr^{\alpha }U(\theta )\leq 1\text{ if }y=1
\end{equation*}%
holds for $k$ small enough such that 
\begin{equation}
kH^{\alpha }U(\theta )\leq 1\text{, where }H=\sqrt{a^{2}+1.}
\end{equation}%
and the proof of the existence of $u_{0,-}(x,y)$, \textit{a }$p_{0}-$\textit{%
subsolution} of problem $P_{a,0,j}$, assumed in Proposition \ref{Propo
Matching} is now completed. 

\noindent With respect to the nondegenerate estimate (\ref%
{Estimate nondegenerate}) of Theorem \ref{Theorem No main}, we point out that the
estimate from above on $\Omega _{+}$ is trivially satisfied since $%
u_{0,-}(x,y)$ is a harmonic function). Thus, the important part is to show
the lower estimate for points $(x,y)\in \Omega _{-}$ which are very near the
part of the boundary in which $y=0.$ For those points, their distance to the
boundary is $\delta (x,y)=y$. Since $r>y,$ from the constructed
subsolution, $u_{0,-}(x,y)=kr^{\alpha }U(\theta ),$ and the fact that we have
\begin{equation*}
c_{1}(\pi -\theta )^{4/3}\leq U(\theta )\leq c_{2}(\pi -\theta )^{4/3},
\end{equation*}%
for $\theta $ sufficiently close to $\pi $, with positive constants $%
c_{1},c_{2}$ (as we can easely prove via the comparison with
sub/supersolutions of the form  $\varphi _{\pm }(x):=A_{\pm }x^{4/3},$ $%
x=\pi -\theta $), we get  the nondegenerate estimate (\ref{Estimate
nondegenerate}) of Theorem \ref{Theorem No main}. Notice that the
assumptions of the uniqueness result Theorem \ref{Thm uniqueness} are
satisfied. $\fin$

\subsection{Construction of a partially flat supersolution}

We make now the structural assumptions, (\ref%
{Hypo j supersolution}), stated in Theorem \ref{TheomMain}. Thus, in particular, we
already know, by Theorems \ref{Theorem No main} and \ref{Thm uniqueness}
that the mere existence of a flat nondegenerate subsolution implies the
uniqueness of the solution. We concentrate our attention now in proving that
the above conditions on $j(x)$ are sufficient conditions for the existence
of a partially flat supersolution, i.e., being flat only on the half of the cathode
region under study $[-a,0]\times \{0\}.$
\bigskip

\noindent \textit{Proof of Theorem 29.} Since we only need to work with the inequality $\geq $ in the equation, it
is enough to use (as before) that 
\begin{equation*}
j(r,\theta )=\frac{A}{(-r\cos \theta )^{\beta }}\geq V_{0}\text{ for }x\in
(-a,0)\text{.}
\end{equation*}%
\noindent We will search for the supersolution by some matching arguments, now in
three different regions, as indicated in the Figure \ref{Fig descomposicion
en 3 zonas} .

\begin{figure}[htp]
\begin{center}
\includegraphics[width=13cm]{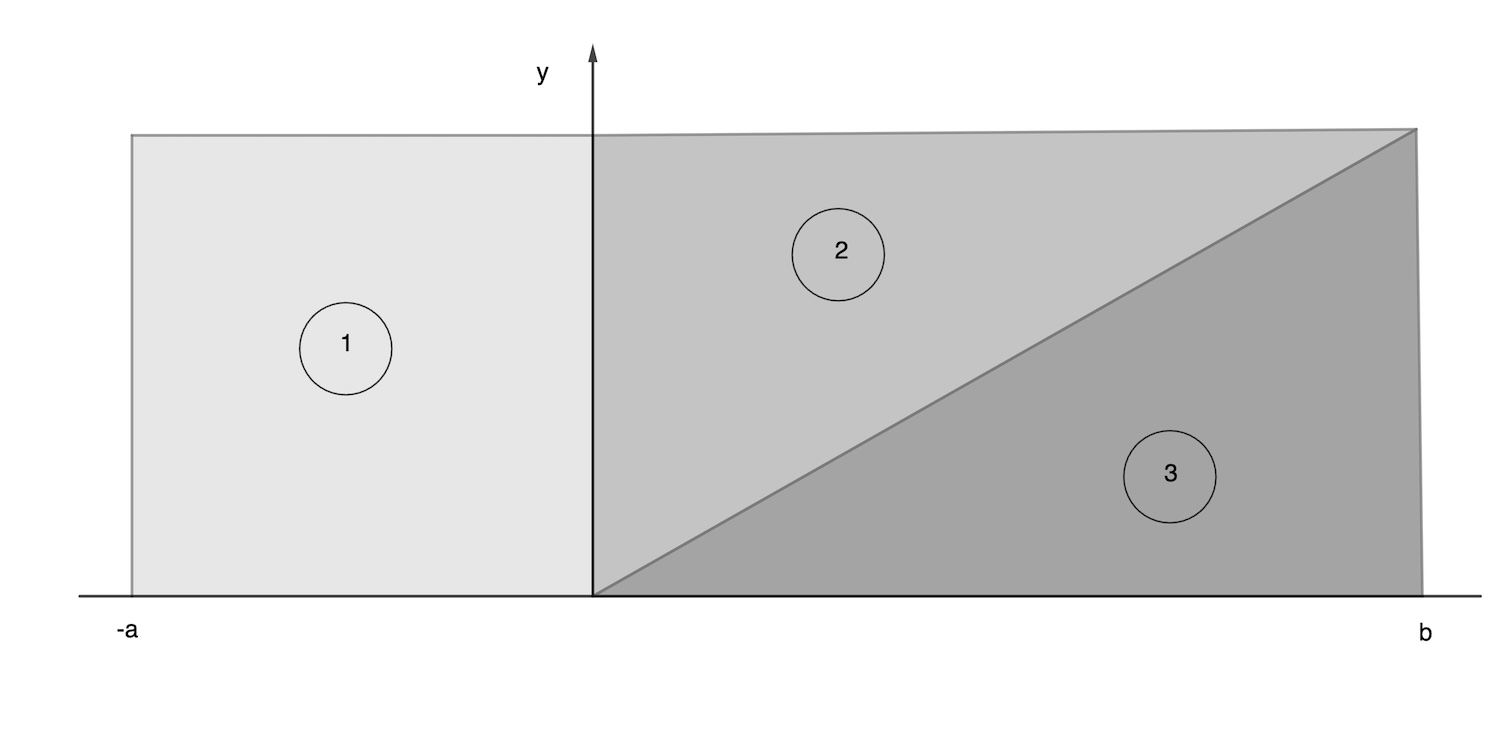}\\ 
\caption{Decomposition in subdomains.}
\label{Fig descomposicion en 3 zonas}
\end{center}
\end{figure}
\noindent In the first region $\Omega _{-}$ we will construct the supersolution $%
u_{-}^{0}(x,y)$ in the same form than the subsolution, i.e.,%
\begin{equation*}
u_{-}^{0}(x,y)=\overline{\phi }(r,\theta )=\overline{k}r^{\overline{\alpha }}%
\overline{U}(\theta )\text{, for some }\overline{k}>0\text{, }\overline{%
\alpha }>1,
\end{equation*}%
and thus with $U(\theta )$ solution of

\begin{equation}
\begin{array}{lr}
-\overline{U}^{\prime \prime }(\theta )+\frac{V_{0}}{\sqrt{(\overline{U}%
(\theta )}}=\lambda ^{\ast }\overline{U}(\theta ) & \theta \in (\frac{\pi }{2%
},\pi ).%
\end{array}
\label{Ec ODE2}
\end{equation}%
So that no boundary layer needs to be taken into account now. This means
that, in terms of the proof of Theorem \ref{Theorem No main}, $\overline{U}%
(\theta )$ is very similar to $v_{2}(\theta ).$ Nevertheless, since we
should ensure that $u_{-}^{0}(x,y)\geq u_{0,-}(x,y)$ we will take the
parameter $\lambda ^{\ast }=(\overline{\alpha })^{2}=\frac{1}{2\overline{R}%
^{2}}\gamma (r_{F})^{2}$, such that $\overline{\alpha }\leq \alpha $ (the
exponent used in the definition of the subsolution), which is obtained if we
take $\frac{\pi }{2}\leq \overline{R}\leq R$ (with $R$ the value choosen for
the part of the subsolution given by $v_{2}(\theta )$). In this way, we will
also have that 
\begin{equation*}
1-\delta \leq \overline{U}(\frac{\pi }{2})\leq 1,\text{ and }-\delta \leq 
\overline{U}^{\prime }(\frac{\pi }{2})\leq 0,
\end{equation*}%
for some $\delta >0$ small enough. Note that even if $V_{0}\leq
V_{\varepsilon }$, the solution structure of these semilinear problems
allows us to take $\overline{k}>k$ large enough (with $k$ the value taken
for the subsolution) in order to have $u_{-}^{0}(x,y)\geq u_{0,-}(x,y)$ on $%
\Omega _{-}$.

\noindent The matching with the second region will follow a different
argument. In fact, on $\Omega _{+}$, the supersolution $u_{+}^{0}$ must be
superharmonic. This will ensure that $u_{+}^{0}(x,y)\geq u_{0,+}(x,y)$ also
on $\Omega _{+}=\Omega _{+}^{2}\cup \Omega _{+}^{3}$ (see Figure \ref{Fig
descomposicion en 3 zonas}), since $u_{0,+}$ was taken as subharmonic and
the inequality $u_{+}^{0}\geq u_{0,+}$ is realized on the boundary of $%
\Omega _{+}$. We start by searching $u_{+}^{0}(x,y)$ on $\Omega _{+}^{2}$ in
the form%
\begin{equation*}
u_{+,2}^{0}(x,y)=Kr^{\alpha }\sin (\alpha \theta ),\text{ }\theta \in
(\theta _{b},\frac{\pi }{2}),
\end{equation*}%
with $\tan \theta _{b}=b$. It is not difficult to check that 
\begin{equation*}
-\Delta u_{+,2}^{0}=K\alpha (\alpha -1)r^{\alpha -2}\sin (\alpha \theta
)\geq 0.
\end{equation*}%
The matching with $u_{-}^{0}(x,y)$ is a delicate question. For a $H^{1}$%
-matching we must have:%
\begin{equation*}
\left \{ 
\begin{array}{cc}
u_{-}^{0}(0,y)=u_{+,2}^{0}(0,y) & y\in (0,1), \\ 
\nabla u_{-}^{0}(0,y)=\nabla u_{+,2}^{0}(0,y) & y\in (0,1).%
\end{array}%
\right. 
\end{equation*}%
The first condition holds if 
\begin{equation*}
\frac{k}{K}=\frac{\sin (\frac{\alpha \pi }{2})}{\overline{U}(\frac{\pi }{2})}%
.
\end{equation*}%
The second condition (once we allow the formation of singularities with a
good sign, i.e., with a non-negative generated implicit distribution) leads
to%
\begin{equation*}
\overline{U}^{\prime }(\frac{\pi }{2})\geq \frac{K}{k}\alpha \cos (\frac{%
\alpha \pi }{2})=\alpha \overline{U}(\frac{\pi }{2})\frac{\cos (\frac{\alpha
\pi }{2})}{\sin (\frac{\alpha \pi }{2})}=\alpha \overline{U}(\frac{\pi }{2}%
)\cot (\frac{\alpha \pi }{2}).
\end{equation*}%
But, $\frac{\alpha \pi }{2}\in (\frac{\pi }{2},\frac{2\pi }{3})\subset
(0,\pi )$ so that $\cot (\frac{\alpha \pi }{2})<0$ (see the Figure \ref{Fig
cotangente}).

\begin{figure}[htp]
\begin{center}
\includegraphics[width=7cm]{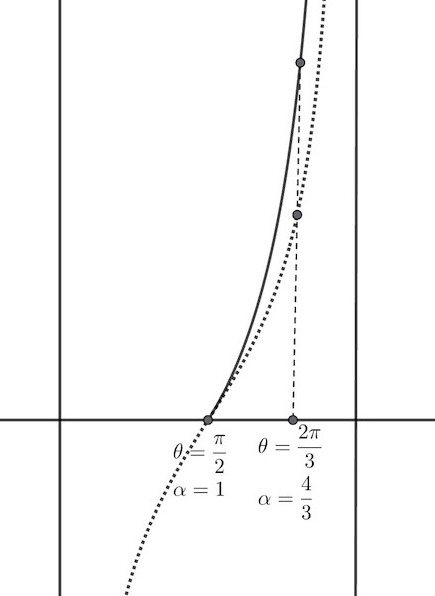}\\[0pt]
\end{center}
\caption{Representation of the functions -$\protect \alpha \cot (\frac{%
\protect \alpha \protect \pi}{2})$}
\label{Fig cotangente}
\end{figure}

\noindent Thus we must investigate for which $\alpha \in (1,\frac{4}{3}]$ we
have some kind of Robin type boundary inequality on $\overline{U}(\theta )$
for $\theta =\frac{\pi }{2}:$ 
\begin{equation}
\frac{-\overline{U}^{\prime }(\frac{\pi }{2})}{\overline{U}(\frac{\pi }{2})}%
\leq \alpha (-\cot (\frac{\alpha \pi }{2})).  \label{Robin inequality}
\end{equation}%
This condition trivially holds if $\alpha =\frac{4}{3}$ (since, in that
case, $\overline{U}^{\prime }(\frac{\pi }{2})=0$). Then by continuity in $%
\alpha $, there exists $\alpha _{0}\in (1,\frac{4}{3})$ such that (\ref%
{Robin inequality}) holds for any $\alpha \in (1,\frac{4}{3})$. This
justifies the assumption (\ref{Hypo j supersolution}) in Theorem \ref%
{Theorem Supersolution} (and in Theorem \ref{TheomMain}).

\noindent Finally, in the third region $\Omega _{+}^{3}$ , we extend the
above function $u_{+,2}^{0}(x,y)$, by a kind of interpolation argument
(remember that we must have the boundary condition $u_{+,3}^{0}(b,y)\geq y$,
for any $y\in (0,1)$). Then we define 
\begin{equation*}
u_{+,3}^{0}(x,y)=C_{3}r\sin (\alpha \theta )+\widehat{C}_{3}r^{\alpha }\sin
(\alpha \theta _{b})\text{ \  \ if }(x,y)\in \Omega _{+,3},
\end{equation*}%
with $\tan \theta _{b}=b.$ It is easy to prove that $u_{+,3}^{0}(x,y)$ is a
super-harmonic function ($-\Delta u_{+,3}^{0}\geq 0$) and that (if $b$ is
large enough) the positive constants $C_{3}$ and $\widehat{C}_{3}$ can be
chosen such that $u_{+,3}^{0}(x,y)$ satisfies the correct inequalities on
the boundaries 
\begin{equation*}
\left \{ 
\begin{array}{cc}
u_{+,3}^{0}(b,y)\geq y, & \text{for any }y\in (0,1), \\[0.15cm]
u_{+,3}^{0}(x,0)\geq 0 & \text{for any }x\in (0,b),%
\end{array}%
\right. 
\end{equation*}%
and it correctly matches with $u_{+,2}^{0}(x,y)$ (generalizing, at most, a
\textquotedblleft good signed\textquotedblright \ measure) on the matching
boundary $\theta =\theta _{b}$. This proves Theorem \ref{Theorem
Supersolution} (and also Theorem \ref{TheomMain}, except the proof of the
uniqueness of solution which will be presented in the next subsection).$%
_{\blacksquare }$

\subsection{Proof of the uniqueness of nondegenerate solutions}

The proof of Theorem \ref{Thm uniqueness} will be a consequence of some
previous results concerning the parabolic problem
\begin{equation}
PP_{a,b,j,u_{0}}=\left\{ 
\begin{array}{cc}
u_{t}-\Delta u+\frac{j(x)}{\sqrt{u}}\chi _{\{u>0\}}=0 & t>0,\text{ }x\in (-a,b),\text{ }%
y\in (0,1), \\ [.3cm]
u(t,x,0)=0 & t>0,\text{ }x\in (-a,b), \\ [.15cm]
u(t,x,1)=1 & t>0,\text{ }x\in (-a,b), \\ [.15cm]
u(t,-a,y)=y^{4/3} & t>0,\text{ }y\in (0,1), \\ [.15cm]
u(t,b,y)=y^{{}} & t>0,\text{ }y\in (0,1), \\ [.15cm]
u(0,x,y)=u_{0}(x,y) & x\in (-a,b),\text{ }y\in (0,1).%
\end{array}%
\right. 
\end{equation}%

Notice that the above parabolic problem is not related to the evolution problem associated with the electron beam formulation. It would be related to some associated evolution problem if the modelling under study were concerned with the Chemical Engineering problem dealing with a chemical kinetics with a reaction of order -1/2.  Here, we are only interested in the fact that this artificial parabolic problem leads us (under suitable assumptions) to the uniqueness of solutions of the mathematical model for a stationary electron beam.

The existence of solutions to problem $PP_{a,b,j,u_{0}}$ is an easy variation of results dealing with parabolic quenching type problems (see, e.g., \cite{Davila-Montenegro},  \cite%
{DaM2003}, \cite{DiGiacomoni} and \cite{Dao-Diaz}). We point out that no
flat condition is required in this formulation. As in the elliptic problem,
we can prove that for certain initial data, the solutions are non-degenerate
in the sense that for any $t>0$ the solution $u$ is such that $u(t)$ belongs
to the following set

\begin{equation*}
\mathcal{M}(\nu ):=\Big\lbrace u\in L^{2}(\Omega ;\delta )\;\big\vert\;\text{%
such\ that }u(x,y)\geq C\delta (x,y)^{\nu }\quad \text{in }\Omega \text{,
for some }C>0\Big\rbrace,\text{ }\nu \in \left( 0,\frac{4}{3}\right] .
\end{equation*}

The following result is an adaptation to this framework of the main theorem of \cite%
{DiGiacomoni}. It gives the continuous dependence of solutions for the initial data (implying, obviously, the uniqueness of solutions),
as well as a smoothing effect concerning the initial datum. In the proof,
we will use some Hilbertian techniques. We consider two initial
data, and we will prove that if the respective solutions are such that $%
u(t),v(t)\in \mathcal{M}(\nu )$ (i.e. with $\delta ^{-\nu }u,\delta ^{-\nu
}v\geq C$)\ then we can estimate the ${L}^{2}(\Omega )$-norm of $\delta
^{-\gamma }[u(t)-v(t)]_{+}$ for some suitable $\gamma \in (0,1]$ in terms of
the $L^{2}(\Omega ;\delta )$-norm of $[u_{0}-v_{0}]_{+}$. Notice that this
automatically implies an estimate on the ${L}^{2}(\Omega )$-norm of $%
[u(t)-v(t)]_{+}$.

\begin{theorem}
\label{TheoremContDependTPatabolic} Let $u_{0},v_{0}\in L^{1}(\Omega )\cap
L^{2}(\Omega ;\delta ).$ Let $u$, $v$ be weak solutions of $PP_{a,b,j,u_{0}}$
and $PP_{a,b,j,v_{0}}$, respectively such that $u(t),v(t)\in \mathcal{M}(\nu
)$ for some $\nu \in \left( 0,\frac{4}{3}\right] $. Then, for any $t\in
(0,\infty )$, we have 
\begin{equation}
\left\Vert \delta ^{-\gamma }[u(t)-v(t)]_{+}\right\Vert _{{L}^{2}(\Omega
)}\leq Ct^{-\frac{2\gamma +1}{4}}\left\Vert [u_{0}-v_{0}]_{+}\right\Vert
_{L^{2}(\Omega ;\delta )},  \label{CDE}
\end{equation}%
with 
\begin{equation}
\gamma :=\min \left\{ \frac{3\nu }{2},1\right\}  \label{gamma}
\end{equation}%
and for some constant $C>0$ independent of $t$. In particular, $u_{0}\leq
v_{0}$ implies that for any $t\in \lbrack 0,+\infty )$, 
\begin{equation*}
u(t,\cdot )\leq v(t,\cdot )\quad \text{ a.e. in }\Omega
\end{equation*}%
and 
\begin{equation}
\left\Vert \delta ^{-\gamma }\left( u(t)-v(t)\right) \right\Vert _{{L}%
^{2}(\Omega )}\leq Ct^{-\frac{2\gamma +1}{4}}\left\Vert
u_{0}-v_{0}\right\Vert _{L^{2}(\Omega ;\delta )}.  \label{CDE2}
\end{equation}
\end{theorem}

\bigskip

Before referring to the proof of Theorem \ref{TheoremContDependTPatabolic},
we will prove that this implies the uniqueness of the positive solution for
the stationary problem $P_{a,b,j}$ presented in Theorem \ref{Thm uniqueness}.

\bigskip

\bigskip

\noindent \textit{Proof of Theorem \ref{Thm uniqueness}.} Let us call $%
u_{\infty }$ and $v_{\infty }$ two possible solutions of $P_{a,b,j}$ in the
class $\mathcal{M}(\nu )$. By taking $u_{0}=u_{\infty }$ and $%
v_{0}=v_{\infty }$ as initial data in $PP_{a,b,j,u_{0}}$ and $%
PP_{a,b,j,v_{0}}$, since $u_{\infty }$ and $v_{\infty }$ are, obviously,
solutions of the respective parabolic problems, we get that $u_{\infty
}-v_{\infty }$ satisfies 
\begin{equation*}
\left\Vert \delta ^{-1}\left( u_{\infty }-v_{\infty }\right) _{+}\right\Vert
_{L^{2}(\Omega )}\leq Ct^{-\frac{2\gamma +1}{4}}\Vert (u_{\infty }-v_{\infty
})_{+}\Vert _{L^{2}(\Omega ,\delta )}.
\end{equation*}%
Making $t\nearrow +\infty $ and reversing the role of $u_{\infty }$ and $%
v_{\infty }$, we get that $u_{\infty }=v_{\infty }.\fin$

\bigskip

The proof of Theorem \ref{TheoremContDependTPatabolic} follows some slight
modifications of the paper \cite{DiGiacomoni}.
Nonetheless, for the sake of completeness, we present here the main lines of
the proof.

Without loss of generality, we can use the notion of \textit{mild solution }%
on $L^{1}(\Omega )$, i.e. $u\in \mathcal{C}([0,T);L^{1}(\Omega )),$ for any $%
T>0$: $j(x)u^{-1/2}\in L^{1}(\Omega \times (0,T))$ and $u$ {fulfills} the
identity 
\begin{equation}
u(\cdot ,t)=S(t)u_{0}(\cdot )-\int_{0}^{t}S(t-s)j(\cdot)(\chi _{\{u>0\}}u^{-1/2}(\cdot ,s))ds,\quad \text{in }L^{1}(\Omega ),  \label{1.2}
\end{equation}%
where $S(t)$ is the $L^{1}(\Omega )$-semigroup corresponding to the Laplace
operator with the corresponding Dirichlet (stationary) boundary conditions
(see, e.g., \cite{Dao-Diaz}). We shall need some well-known auxiliary
results. The first one is a singular version of the Gronwall's inequality
(due to Brezis and Cazenave \cite{Brezis-Cazenave}), which is especially
useful in the study of non-global Lipschitz perturbations of the heat
equation.

\begin{lemma}
\label{GronLemma} Let $T>0$, $A\geq 0$, $0\leq a,b\leq 1$ and let $f$ be a
non-negative function with $f\in L^{p}(0,T)$ for some $p>1$ such that $%
\max \{a,b\}<1/p^{\prime }$ (where $\frac{1}{p}+\frac{1}{p^{\prime }}=1$).
Consider a non-negative function $\varphi \in L^{\infty }(0,T)$ such that,
for almost every $t\in (0,T)$, 
\begin{equation}
\varphi (t)\leq At^{-a}+\int_{0}^{t}(t-s)^{-b}f(s)\varphi (s)\,ds.
\label{Gron1}
\end{equation}%
Then, there exists $C>0$ only depending on $T,a,b,p$ and $\Vert f\Vert
_{L^{p}(0,T)}$ such that, for almost every $t\in (0,T)$, 
\begin{equation}
\varphi (t)\leq ACt^{-a}.  \label{Gron2}
\end{equation}
\end{lemma}

\bigskip

We shall also use some regularizing effects properties satisfied by the
semigroup $S(t)$ of the heat equation with zero Dirichlet boundary conditions (since we will apply it to the difference of two solutions of the parabolic problem $PP_{a,b,j,u_{0}}$)
due to different authors, among them \cite{Veron}, \cite{Davila-Montenegro} and \cite{Souplet}: see details, in \cite{DiGiacomoni}.

\begin{lemma}
\label{SGL} ~~

\begin{enumerate}
\item \label{SGI1} There exists $C>0$ such that, for any $t>0$ and any $%
u_{0}\in L^{2}(\Omega )$, 
\begin{equation}
\Vert \nabla S(t)u_{0}\Vert _{L^{2}(\Omega )}\leq Ct^{-\frac{1}{2}}\Vert
u_{0}\Vert _{L^{2}(\Omega )}.  \label{SGP1}
\end{equation}

\item \label{SGI2} There exists $C>0$ such that, for any $t>0$ and any $%
u_{0}\in L^{1}(\Omega )$, 
\begin{equation}
\Vert S(t)u_{0}\Vert _{L^{2}(\Omega )}\leq Ct^{-\frac{N}{4}}\Vert u_{0}\Vert
_{L^{1}(\Omega )}.  \label{SGP2}
\end{equation}

\item \label{SGI3} There exists $C>0$ such that, for any $t>0$, any $m\in
(0,1]$ and any $u_{0}\in L^{2}(\Omega ;\delta ^{2m})$, 
\begin{equation}
\Vert S(t)u_{0}\Vert _{L^{2}(\Omega )}\leq Ct^{-\frac{m}{2}}\Vert u_{0}\Vert
_{L^{2}(\Omega ,\delta ^{2m})}.  \label{SGP3}
\end{equation}

\item \label{SGI4} There exists $C>0$ such that, for any $t>0$, any $p\in
\lbrack 1,+\infty )$ and any $u_{0}\in L^{p}(\Omega ,\delta )$, 
\begin{equation}
\Vert S(t)u_{0}\Vert _{L^{p}(\Omega )}\leq Ct^{-\frac{1}{2p}}\Vert
u_{0}\Vert _{L^{p}(\Omega ,\delta )}.  \label{SGP4}
\end{equation}
\end{enumerate}
\end{lemma}

Finally, we can give the main arguments of the proof.

\bigskip

\noindent \noindent \textit{Proof of Theorem \ref%
{TheoremContDependTPatabolic}. }By the constant variations formula, we know
that for any $t\in \lbrack 0,T]$, 
\begin{equation}
u(t)-v(t)=S(t)(u_{0}-v_{0})+\int_{0}^{t}S(t-s)\left( h(u(s))-h(v(s))\right)
\,ds\quad \text{ in }\Omega ,
\end{equation}%
where $h(x,u):=j(x)u^{-1/2}$. By the convexity of the function $u\mapsto
u^{-1/2}$ and the assumption that $u(t),v(t)\in \mathcal{M}(\nu )$, we
deduce that 
\begin{equation}
h(x,u)-h(x,v)\leq Cj(x)\delta ^{-3\nu /2}(u-v)_{+}\quad \text{ in }\Omega .
\end{equation}%
Thus, if we denote $w:=u-v$, we get for any $\tau ,t\in \lbrack 0,T]$ with $%
\tau \leq t$ 
\begin{equation}
w_{+}(t)\leq S(t-\tau )w_{+}(\tau )+C\int_{\tau }^{t}S(t-s)(j(x)\delta
^{-3\nu /2}w_{+}(s))\,ds.  \label{VC1}
\end{equation}%
We multiply \eqref{VC1} by the weight $\delta ^{-\gamma }$, with $\gamma \in
\lbrack 0,1]$ to be chosen later, and take the $L^{2}$-norms. Then, 
\begin{equation*}
\Vert \delta ^{-\gamma }w_{+}(t)\Vert _{L^{2}(\Omega )}\leq \Vert \delta
^{-\gamma }S(t-\tau )w_{+}(\tau )\Vert _{L^{2}(\Omega )}+C\int_{\tau
}^{t}\Vert S(t-s)j(x)\delta ^{-[(\beta +1)\nu +\gamma ]}w_{+}(s)\Vert
_{L^{2}(\Omega )}\,ds.
\end{equation*}%
Let us fix $s,t>0$ and let us call $\psi :=S(t-s)j(x)\delta ^{-(\beta +1)\nu
}w_{+}(s).$ Then, by H\"{o}lder inequality, 
\begin{equation*}
\Vert \delta ^{-\gamma }\psi \Vert _{L^{2}(\Omega )}^{2}=\int_{\Omega }\frac{%
\psi ^{2}}{\delta ^{2\gamma }}\,dx\leq \left( \int_{\Omega }\frac{\psi ^{2}}{%
\delta ^{2}}\,dx\right) ^{\gamma }\left( \int_{\Omega }\psi ^{2}\,dx\right)
^{1-\gamma }
\end{equation*}%
(note that the limit cases $\gamma \equiv 0$ and $\gamma \equiv 1$ are
allowed). Then, applying the Hardy inequality, 
\begin{equation*}
\Vert \delta ^{-\gamma }\psi \Vert _{L^{2}(\Omega )}\leq C\Vert \nabla \psi
\Vert _{L^{2}(\Omega )}^{\gamma }\Vert \psi \Vert _{L^{2}(\Omega
)}^{1-\gamma }.
\end{equation*}%
By property \ref{SGI1} of Lemma \ref{SGL}, for $\frac{t-s}{2}$, we get 
\begin{equation}
\Vert \delta ^{-[(\beta +1)\nu +\gamma ]}S(t-s)w_{+}(s)\Vert _{L^{2}(\Omega
)}\leq C(t-s)^{-\frac{\gamma }{2}}\Vert S\left( \frac{t-s}{2}\right)
j(x)\delta ^{-(\beta +1)\nu }w_{+}(s)\Vert _{L^{2}(\Omega )}.
\end{equation}%
Analogously, using property \ref{SGI4} of Lemma \ref{SGL}, and that $j\in
L^{2}(\Omega )$ (remember the assumptions made in Theorem \ref{Theorem No
main}) 
\begin{equation}
\Vert \delta ^{-\gamma }S(t)w_{+}(0)\Vert _{L^{2}(\Omega )}\leq Ct^{-\frac{%
\gamma }{2}}\Vert S\left( \frac{t}{2}\right) j(x)w_{+}(0)\Vert
_{L^{2}(\Omega )}\leq Ct^{-\left( \frac{\gamma }{2}+\frac{1}{4}\right)
}\Vert w_{+}(0)\Vert _{L^{2}(\Omega ;\delta )}.
\end{equation}%
In order to apply the singular Gronwall's inequality, we must relate the
weights $\delta ^{-\gamma }$ and $\delta ^{-3\nu /2}$ keeping in mind that $%
\gamma \in \lbrack 0,1]$. To do that, we apply property \ref{SGI3} of Lemma %
\ref{SGL} for some $m\in \lbrack 0,1]$. We shall take 
\begin{equation}
3\nu /2=\gamma +m.  \label{choice of m}
\end{equation}%
Indeed, if $(\beta +1)\nu \in (1,2]$, then we take $\gamma =1$, $m=3\nu /2-1$
and we apply point \ref{SGI3} of Lemma \ref{SGL} to the initial datum: 
\begin{equation}
\Vert S\left( \frac{t-s}{2}\right) j(x)\delta ^{-(\beta +1)\nu
}w_{+}(s)\Vert _{L^{2}(\Omega )}=\Vert S\left( \frac{t-s}{2}\right)
j(x)\delta ^{-(m+1)}w_{+}(s)\Vert _{L^{2}(\Omega )}\leq C(t-s)^{-\frac{m}{2}%
}\Vert \delta ^{-\gamma }w_{+}(s)\Vert _{L^{2}(\Omega )}.
\end{equation}%
On the other hand, if $3\nu /2\in \lbrack 0,1]$, we can take $\gamma =3\nu
/2 $ and thus, since $S(t-s)$ is a contraction in $L^{2}(\Omega )$, we get 
\begin{equation}
\Vert S\left( \frac{t-s}{2}\right) j(x)\delta ^{-(\beta +1)\nu
}w_{+}(s)\Vert _{L^{2}(\Omega )}=\Vert S\left( \frac{t-s}{2}\right)
j(x)\delta ^{-\gamma }w_{+}(s)\Vert _{L^{2}(\Omega )}\leq \Vert \delta
^{-\gamma }w_{+}(s)\Vert _{L^{2}(\Omega )},
\end{equation}%
which corresponds to \eqref{choice of m} with $m=0$. In other words, 
\begin{equation*}
\gamma =\min \{{1,3\nu /2}\}
\end{equation*}%
and 
\begin{equation*}
m=\max \{3\nu /2-1,0\}.
\end{equation*}%
Collecting the previous inequalities, we arrive to 
\begin{equation*}
\Vert \delta ^{-\gamma }w_{+}(t)\Vert _{L^{2}(\Omega )}\leq Ct^{-\frac{%
2\gamma +1}{4}}\Vert w_{+}(0)\Vert _{L^{2}(\Omega ;\delta
)}+C\int_{0}^{t}(t-s)^{-\frac{m}{2}}\Vert \delta ^{-\gamma }w_{+}(s)\Vert
_{L^{2}(\Omega )}.
\end{equation*}%
Thus, we can apply Lemma \ref{GronLemma} with $a=\frac{2\gamma +1}{4}\in %
\left[ \frac{1}{4},\frac{3}{4}\right] $, $b=\frac{m}{2}$ and $A=C\Vert
w_{+}(0)\Vert _{L^{2}(\Omega ;\delta )}$ to deduce that 
\begin{equation*}
\Vert \delta ^{-\gamma }w_{+}(t)\Vert _{L^{2}(\Omega )}\leq Ct^{-\frac{%
2\gamma +1}{4}}\Vert w_{+}(0)\Vert _{L^{2}(\Omega ;\delta )}.
\end{equation*}
\fineq

\begin{remark}
As in the proof of Theorem \ref{Thm uniqueness}, if $u_{\infty }$ is the non-degenerate solution of
problem $P_{a,b,j}$, for any $u_{0}\in L^{2}(\Omega ;\delta )$ leading to a
non-degenerate solution $u(t)$ of the parabolic problem $PP_{a,b,j,u_{0}}$, we get a
quantitative estimate on the asymptotic stability of $u_{\infty }$ in the
class of non-degenerate stationary solutions, since, for any $t\in (0,\infty
)$, we have 
\begin{equation*}
\left\Vert \delta ^{-\gamma }[u(t)-u_{\infty }]\right\Vert _{{L}^{2}(\Omega
)}\leq Ct^{-\frac{2\gamma +1}{4}}\left\Vert u_{0}-u_{\infty }\right\Vert
_{L^{2}(\Omega ;\delta )},
\end{equation*}%
with $\gamma :=\min \left\{ \frac{3\nu }{2},1\right\} $, for some
constant $C>0$ independent of $t$.
\end{remark}

\par
\noindent 
{\bf Acknowledgement}
Research partially supported by the project PID2023-146754NB-I00 funded by
MCIU/AEI/10.13039/501100011033 and FEDER, EU.
MCIU/AEI/10.13039/-501100011033/FEDER, EU. The author thanks Ha\"{\i}m Brezis and
Joel Lebowitz for several conversations on this subject. Thanks are also extended to my brother, Gregorio Díaz, and my friend, Laurent Véron, for their careful reading of an earlier version of the manuscript. They suggested several improvements to the exposition. In particular, L. Véron showed me an alternative proof of Theorem 2.

\begin{figure}[htp]
\begin{center}
\includegraphics[width=9cm]{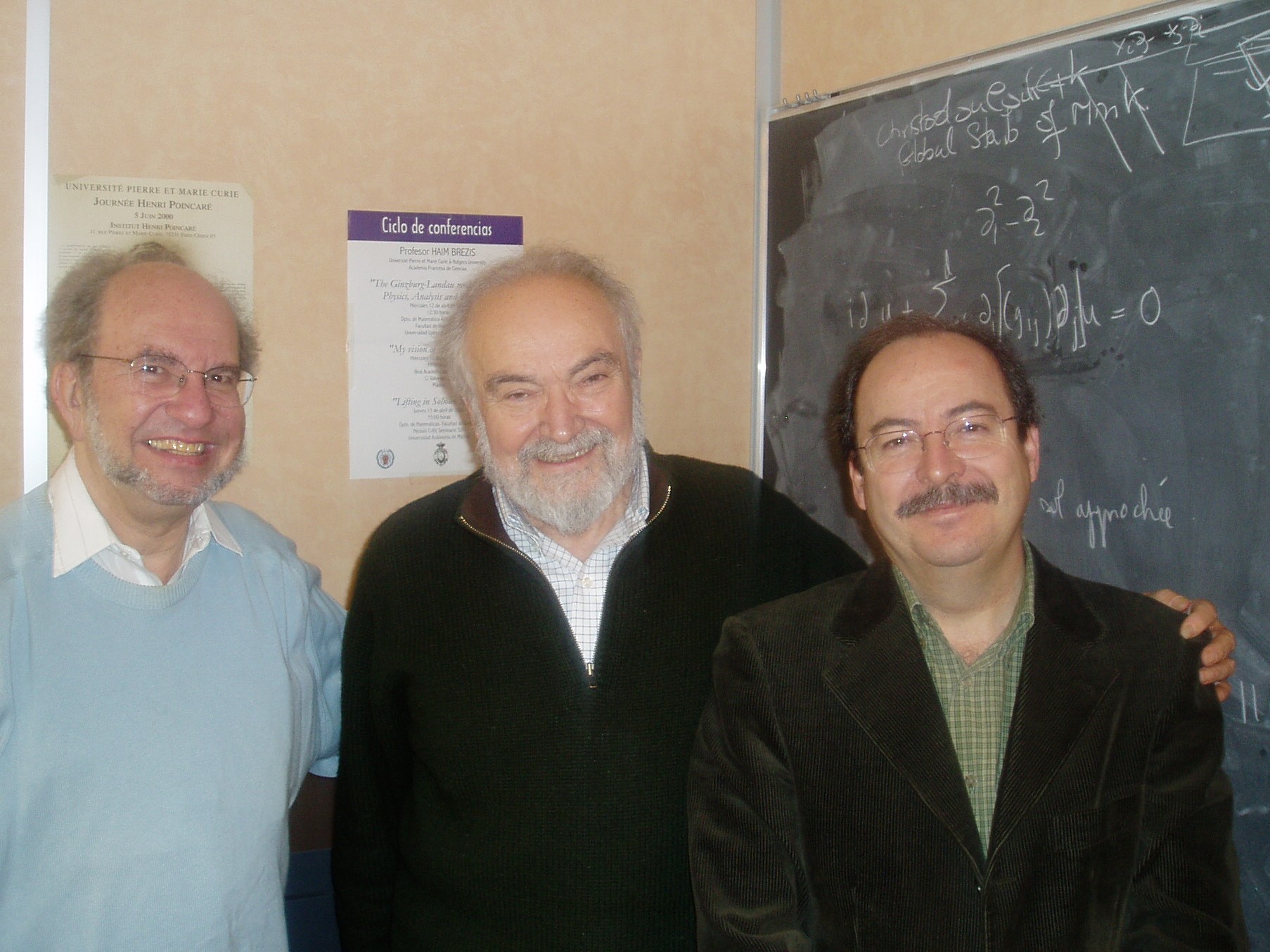}\\ 
\caption{Haïm Brezis, Joel Lebowitz and the author, in Paris, 2005.}
\end{center}
\end{figure}

\bibliographystyle{amsplain}
\bibliography{CRASBrezis2025}

%
%
%

\end{document}